\documentclass[11pt]{article}
\usepackage{amsfonts,amsbsy,amsmath,amssymb}
\usepackage{graphics}
\usepackage[lmargin=3.0cm, rmargin=3.0cm, tmargin=3.0cm, bmargin=3.0cm]{geometry}

\numberwithin{equation}{section}

\newcommand{\iso}{\cong} 

\newcommand{\scaps}[1]{{\scshape #1}}
\newcommand{\bfm}{\textbf}
\newcommand{\mcal}{\mathcal}

\newcommand{\mbs}{\boldsymbol}

\newcommand{\A}{\mcal{A}}
\newcommand{\B}{\mcal{B}}
\newcommand{\C}{\mcal{C}}
\newcommand{\D}{\mcal{D}}

\renewcommand{\H}{\mcal{H}}

\renewcommand{\L}{\mcal{L}}

\renewcommand{\P}{\mcal{P}}

\newcommand{\X}{\mcal{X}}
\newcommand{\Y}{\mcal{Y}}

\newcommand{\bnd}[2]{{\rm bnd}_{#2}#1}
\newcommand{\half}{\mbox{\small $\frac{1}{2}$}}
\newcommand{\sm}{\setminus}

\def\proof{\par\noindent{{\sl Proof.\ }}}
\def\QED{$\blacksquare$}
\def\inQED{$\square$}

\newenvironment{statement}
{
\refstepcounter{equation} 
\ \\
\noindent
\begin{it}
\noindent
\emph{\bfm{~\theequation.}}
}
{\end{it}\ \\}

\newcommand{\sref}[1]{\bfm{{\ref{#1}}}}

\makeatletter
\def\section{\@ifstar\unnumberedsection\numberedsection}
\def\numberedsection{\@ifnextchar[
  \numberedsectionwithtwoarguments\numberedsectionwithoneargument}
\def\unnumberedsection{\@ifnextchar[
  \unnumberedsectionwithtwoarguments\unnumberedsectionwithoneargument}
\def\numberedsectionwithoneargument#1{\numberedsectionwithtwoarguments[#1]{#1}}
\def\unnumberedsectionwithoneargument#1{\unnumberedsectionwithtwoarguments[#1]{#1}}
\def\numberedsectionwithtwoarguments[#1]#2{%
  \ifhmode\par\fi
  \removelastskip
  \vskip 3ex\goodbreak
  \refstepcounter{section}%
  \noindent
  \leavevmode
  \begingroup
  \bfseries
  \S \thesection\ 
  #2.\quad
  \endgroup
  \addcontentsline{toc}{section}{%
    \protect\numberline{\thesection}%
    #1}%
  }
\def\unnumberedsectionwithtwoarguments[#1]#2{%
  \ifhmode\par\fi
  \removelastskip
  \vskip 3ex\goodbreak
  \noindent
  \leavevmode
  \begingroup
  \bfseries
  #2.\quad
  \endgroup
  \addcontentsline{toc}{section}{%
    #1}%
  }

\makeatletter
\def\subsection{\@ifstar\unnumberedsubsection\numberedsubsection}
\def\numberedsubsection{\@ifnextchar[
  \numberedsubsectionwithtwoarguments\numberedsubsectionwithoneargument}
\def\unnumberedsubsection{\@ifnextchar[
  \unnumberedsubsectionwithtwoarguments\unnumberedsubsectionwithoneargument}
\def\numberedsubsectionwithoneargument#1{\numberedsubsectionwithtwoarguments[#1]{#1}}
\def\unnumberedsubsectionwithoneargument#1{\unnumberedsubsectionwithtwoarguments[#1]{#1}}
\def\numberedsubsectionwithtwoarguments[#1]#2{%
  \ifhmode\par\fi
  \removelastskip
  \vskip 3ex\goodbreak
  \refstepcounter{subsection}%
  \noindent
  \leavevmode
  \begingroup
  \bfseries
  \S \thesubsection\ 
  #2.\quad
  \endgroup
  \addcontentsline{toc}{subsection}{%
    \protect\numberline{\thesubsection}%
    #1}%
  }
\def\unnumberedsubsectionwithtwoarguments[#1]#2{%
  \ifhmode\par\fi
  \removelastskip
  \vskip 3ex\goodbreak
  \noindent
  \leavevmode
  \begingroup
  \bfseries
  #2.\quad
  \endgroup
  \addcontentsline{toc}{subsection}{%
    #1}%
  }
\makeatother

\makeatletter
\def\subsubsection{\@ifstar\unnumberedsubsubsection\numberedsubsubsection}
\def\numberedsubsubsection{\@ifnextchar[
  \numberedsubsubsectionwithtwoarguments\numberedsubsubsectionwithoneargument}
\def\unnumberedsubsubsection{\@ifnextchar[
  \unnumberedsubsubsectionwithtwoarguments\unnumberedsubsubsectionwithoneargument}
\def\numberedsubsubsectionwithoneargument#1{\numberedsubsubsectionwithtwoarguments[#1]{#1}}
\def\unnumberedsubsubsectionwithoneargument#1{\unnumberedsubsubsectionwithtwoarguments[#1]{#1}}
\def\numberedsubsubsectionwithtwoarguments[#1]#2{%
  \ifhmode\par\fi
  \removelastskip
  \vskip 3ex\goodbreak
  \refstepcounter{subsubsection}%
  \noindent
  \leavevmode
  \begingroup
  \bfseries
  \S \thesubsubsection\ 
  #2.\quad
  \endgroup
  \addcontentsline{toc}{subsubsection}{%
    \protect\numberline{\thesubsubsection}%
    #1}%
  }
\def\unnumberedsubsubsectionwithtwoarguments[#1]#2{%
  \ifhmode\par\fi
  \removelastskip
  \vskip 3ex\goodbreak
  \noindent
  \leavevmode
  \begingroup
  \bfseries
  #2.\quad
  \endgroup
  \addcontentsline{toc}{subsubsection}{%
    #1}%
  }
\makeatother

\begin{document}
\begin{center}
{\Large \scaps{Subdivisions in apex graphs}}\\ \ \\
Elad Aigner-Horev\\
elad.horev@uni-hamburg.de\\ \ \\
{\small 
Department of Mathematics\\
Hamburg university}
\end{center}

{\small
\noindent
\emph{\bfm{Abstract.}} The Kelmans-Seymour conjecture states that the $5$-connected nonplanar graphs contain a subdivided $K_{_5}$. 
Certain questions of Mader propose a ``plan'' towards a possible resolution of this conjecture. One part of this plan is to show that a $5$-connected nonplanar graph containing $K^-_{_4}$ or $K_{_{2,3}}$ as a subgraph has a subdivided $K_{_5}$. Recently, Ma and Yu showed that  a $5$-connected nonplanar graph containing $K^-_{_4}$ as a subgraph has a subdivided $K_{_5}$. We take interest in $K_{_{2,3}}$ and prove that a $5$-connected nonplanar apex graph containing $K_{_{2,3}}$ as a subgraph has a subdivided $K_{_5}$.\\ 


\noindent 
\scaps{Preamble.} Whenever possible notation and terminology are that of~\cite{diestel}.
Throughout, a graph is always simple, undirected, and finite. $G$ always denotes a graph. A subdivided $G$ is denoted $TG$ and the \emph{branch} vertices of a $TG$ are its vertices at least $3$-valent. A $TG$ is said to be \emph{rooted} at a prescribed vertex set $X$, if the latter coincides with its branch vertices.} 

\section{Introduction}
A refinement of Kuratowski's theorem postulated by the Kelmans-Seymour conjecture ($1975$) is that: {\sl the $5$-connected nonplanar graphs contain a $TK_{_5}$}.  Certain questions of Mader~\cite{Mader1} seem to suggest a two phase plan towards a possible resolution of the Kelmans-Seymour conjecture.\\

\noindent
\scaps{Mader's ``Plan''.} \\
\emph{Phase I.} \emph{Find a small set of graphs $\H$ such that every graph with at least $\left\lfloor \frac{5}{2} n \right\rfloor$ edges contains a $TK_{_5}$ or a member of $\H$ as a subgraph.}\\
\emph{Phase II.} \emph{Prove: for every $H \in \H$, if $G$ is a $5$-connected nonplanar graph with $H$ as a subgraph, then $TK_{_5} \subseteq G$.}\\

The questions of Mader~\cite{Mader1} and the discussion thereafter suggest to examine $\H = \{K^-_{_4}, K_{_{2,3}}\}$, where $K^-_{_4}$ denotes $K_{_4}$ with a single edge removed. 
We are not aware of any mathematics that would indicate as to which phase of the above plan one should consider first. Nevertheless, recently, Ma and Yu~\cite{MaYu1,MaYu2} proved the following.

\begin{statement}\label{MaYuk4}
\bfm{(Ma-Yu~\cite{MaYu1,MaYu2})}\\ 
A $5$-connected nonplanar graph containing $K^-_{_4}$ as a subgraph, contains a $TK_{_5}$. 
\end{statement}

\noindent 
We take interest in $K_{_{2,3}}$ and ask {\sl is it true that a nonplanar $5$-connected graph containing a $K_{_{2,3}}$ as a subgraph contains a $TK_{_5}$?} Using \sref{MaYuk4}, we answer this question in the affirmative for apex graphs, where a graph is apex if it contains a vertex removal of which results in a planar graph. 

\begin{statement}\label{k23}
A $5$-connected nonplanar apex graph containing $K_{_{2,3}}$ as a subgraph contains a $TK_{_5}$. 
\end{statement}

The restriction of \sref{k23} to apex graphs is not artificial. 
Indeed, the conjecture being open for so many years suggests that it does not stand to reason to consider special cases of it. We contend that apex graphs rise quite naturally; for instance, when trying to determine properties of a minimal counter example to the conjecture. 
  
%
%

Our approach to prove \sref{k23} is structural and as a result most of its significant parts do not use the assumption that the graph is apex. The generalization of the parts of the proof that do relay on this assumption coincide with a well-known problem proposed by Yu~\cite{Yu}. To emphasize this, the proof is written in a modular manner 
so that the dependency on planarity (i.e., the apex assumption) can be easily singled out. Next, we make this more precise as we outline the highlights of our approach. \\

\noindent   
\bfm{I.} Essentially, to prove \sref{k23}, we complete a $K_{_{2,3}}$ into a $TK_{_5}$. Such an approach
brings forth a problem which we have come to call the \emph{rooted nonseparating path problem}: {\sl for integer $k>0$, $x,y \in V(G)$, and $Y \subset V(G) \sm \{x,y\}$ (the roots), is there an induced $xy$-path $P$ in $G$ satisfying $\kappa(G-P) \geq k$ and $V(P) \cap Y = \emptyset$?} This problem is at the heart of the Ma-Yu~\cite{MaYu1,MaYu2} proof for \sref{MaYuk4}; it will also be at the heart of our proof for \sref{k23}. In particular, we prove \sref{rooted-nonseparating-path} which we propose to be of independent interest. \\

\noindent
\bfm{II.} The assertion of \sref{rooted-nonseparating-path} allows us to assume that if $K_{_{2,3}} \iso K \subseteq G$, then, essentially,  $G$ is comprised of two parts, i.e., $G = H \cup P$, where $P$ is an induced path linking
the $3$-valent vertices of $K$ meeting no $2$-valent vertex of $K$, and $H =G-P$ is a chain of blocks (i.e., maximal $2$-connected components of $G-P$). The following is then one of our main results. 

\begin{statement}\label{2con-main}
Suppose $G$ is a $5$-connected graph such that 
$K_{_{2,3}} \iso K \subseteq G$. If $G$ has an induced path $P$ linking
the $3$-valent vertices of $K$ meeting no $2$-valent vertex of $K$ such that 
$\kappa(G-P) \geq 2$, then $TK_{_5} \subseteq G$.   
\end{statement}

Note that \sref{2con-main} is not limited to apex nor to nonplanar graphs and in our opinion would serve as a good stepping stone towards generalizing \sref{k23} for nonapex graphs.
The characterization of $2$-connected graphs containing no circuit through $3$ prescribed elements (edges or vertices) is the main engine in our proof of \sref{2con-main} (see \S\ref{frames}). \\

\noindent
\bfm{III.} Sadly, \sref{2con-main} will not be sufficient for us as we do not know how to guarantee existence of such a path. When confined to apex graphs a new structural tool is made available; such is the characterization, obtained by Yu~\cite{Yu}, of pairs $(G,W)$ such that $G$ is a $4$-connected plane graph, $W \subseteq V(G),|W|=4$, and $G$ has no $TK_{_4}$ rooted at $W$ (see \S\ref{salami-and-pie} for details). 

At this point of the proof, we have two descriptions of the graph $G$. One description is that which emerges from applying Yu's result in~\cite{Yu} in certain ways; the second is that which emerges from \sref{2con-main} and \sref{rooted-nonseparating-path} that $G = H \cup P$ where $H$ and $P$ are as above and $\kappa(H) =1$.  A comparison between these two descriptions is then drawn as to arrive at a conclusion of our proof.

Appealing to the result of Yu~\cite{Yu} (and the analysis that it allows us to carry) is essentially the sole place of the proof that relies on planarity and the graph being apex. In~\cite{Yu}, Yu proposes the problem of generalizing his result to pairs $(G,W)$ as above where $G$ is not necessarily planar. 
   
It is worth to note that using this characterization, Yu~\cite{Yu} proved that every $6$-connected apex graph contains a $TK_{_5}$. This was then superseded by Mader~\cite{Mader1} who proceeded to prove that every $6$-connected graph contains a $TK_{_5}$.\\

\noindent
\bfm{IV.} Finally, \sref{MaYuk4} and a variant of the following recent result of Ma and Yu~\cite{MaYu1} will be used extensively throughout our arguments for \sref{k23}. Indeed, these two theorems of Ma and Yu form a useful tool in our analysis. 

\begin{statement}\label{MaYuSep}
\bfm{(Ma-Yu~\cite{MaYu1})}\\ 
Let $G$ be a $5$-connected nonplanar graph such that $G =G_{_1} \cup G_{_2}$, $|V(G_{_1} \cap G_{_2})| = 5$. If $|V(G)| > |V(G_{_2})| \geq 7$ and $G_{_2}$ has a planar representation with $V(G_{_1} \cap G_{_2})$ incident with a common face, then $TK_{_5} \subseteq G$.
\end{statement}

\noindent
\scaps{Organization.} Additional notation, terminology, and required previous results are detailed in \S\ref{pre}. A proof of \sref{2con-main} is provided in \S\ref{proof-2con-main}.
A proof of \sref{k23} is provided in \S\ref{apex-main}.

\section{Preliminaries}\label{pre} Throughout, whenever a claim regarding plane graphs is applied to a planar graph $G$, we always assume that the application is done with respect to some embedding of $G$.\\

\noindent
\scaps{Agreement.} We always assume that $G+e =G$ if $e \in E(G)$. In addition to that, if $e\in E(G)$ and we write $G'=G+e$ and then $G''=G'-e$ in an argument, then the agreement is that $G' = G = G''$. We use this agreement in \S\ref{proof-2con-main} and \S\ref{proof-2con2}. \\

\noindent
\scaps{Subgraphs.} Let $H$ be a subgraph of $G$, denoted $H \subseteq G$.
The boundary of $H$, denoted by $\bnd{H}{G}$ (or simply by $\bnd{H}{}$), is the set of vertices of $H$ incident with $E(G)\sm E(H)$. By $int_G H$ (or simply $int H$) we denote the subgraph induced by $V(H) \sm bnd H$. If $v \in V(G)$, then $N_H(v)$ denotes $N_G(v) \cap V(H)$.  

Given $U \subseteq V(G) \cup E(G)$, we write $U \cap H$ to denote $U \cap (V(H) \cup E(H))$.
We also write $x \in G$ to denote that $x$ is an element of $G$, i.e., an edge or a vertex.\\

\noindent
\scaps{Paths and circuits.} 
For $X,Y \subseteq V(G)$, an $(X,Y)$-path is a simple path with one end in $X$ and the other in $Y$ internally-disjoint of $X \cup Y$. If $X =\{x\}$, we write $(x,Y)$-path.
If $|X|=|Y|=k \geq 1$, then a set of $k$ pairwise disjoint $(X,Y)$-paths is called 
an $(X,Y)$-$k$-\emph{linkage}. Also, if $x \in V(G)$ and $Y \subseteq V(G) \sm \{x\}$, then by $(x,Y)$-$k$-\emph{fan} we mean 
a set of $k \geq 1$ $(x,Y)$-paths with only $x$ as a common vertex.

The \emph{interior} of an $xy$-path $P$ is the set $V(P)\sm\{x,y\}$ and is denoted $int P$. 
For $u,v \in V(P)$, we write $uPv$ to denote the vertex set of the $uv$-subpath of $P$ and is called the $uv$-\emph{segment} of $P$. We write $(uPv)$ to denote $uPv \sm \{u,v\}$, and in a similar manner the segment $(uPv$ and the segment $uPv)$ are defined. 

A $uv$-segment of $P$ \emph{nests} in a $u'v'$-segment of $P$ if $uPv \subseteq u'Pv'$. These segments are said to \emph{overlap} if $x,u,u',v,v',y$ appear in this order along $P$. 

Let $H_{_1} \subseteq G-P$. A vertex $u \in V(P)$ adjacent to $H_{_1}$ is called an \emph{attachment vertex} of $H_{_1}$ on $P$. If $H_{_1}$ has $\geq 2$ attachments on $P$, then the attachments $v_{_1}$, $v_{_2}$ such that $v_{_1}P v_{_2}$ is maximal are called the \emph{extremal} attachments of $H_{_1}$ on $P$; the set $v_{_1}P v_{_2}$ is then called the \emph{extremal segment} of $H_{_1}$ on $P$. 

We say that $H_{_1}$ \emph{nests} in $H_{_2} \subseteq G-P$ with respect to $P$ if the extremal segment of $H_{_1}$ on $P$ nests in the extremal segment of $H_{_2}$ on $P$. These subgraphs are said to \emph{overlap} with respect to $P$ if their extremal segments on $P$ overlap. 

Given a set $A\subseteq V(G) \cup E(G)$, we refer to a circuit containing $A$ as an $A$-\emph{circuit}. For a circuit $C$ embedded in the plane we denote by $int C$ and $ext C$ the regions of the plane bounded by $C$ and appearing in its interior and exterior, respectively. \\

\noindent
\scaps{Bridges.} Let $H \subseteq G$. By $H$-\emph{bridge} we mean either an edge $uv \notin E(H)$ and $u,v \in V(H)$ or a connected component of $G-H$. In the latter case, the $H$-bridge is called \emph{nontrivial}. The vertices of $H$ adjacent to an $H$-bridge $B$ are called the \emph{attachment vertices} of $B$. A $uv$-path internally-disjoint of $H$ with $u,v \in V(H)$, is called an $H$-\emph{ear}.\\ 

\noindent
\scaps{Hammocks and Disconnectors.} A $k$-\emph{hammock} in $G$ is a connected subgraph whose boundary consists of $k \geq 1$ vertices; called {\em ends}. A hammock coinciding with its boundary is called {\sl trivial}. In this sense, a trivial $2$-hammock is an edge. 
For $2$-hammocks we also refer to $G-e$ as a trivial $2$-hammock, where $e \in E(G)$.
Clearly, a graph with no nontrivial $2$-hammocks is $3$-connected, and so on.
   
Let $G$ have $\kappa(G) =2$. Then $G$ is a union of $2$-hammocks.  
An end-to-end path in a $2$-hammock will be called a \emph{through path}. Being not 
necessarily $2$-connected, a $2$-hammock becomes such after an ear linking its ends is appended. Equivalently, 

\begin{equation}\label{through-path}
\mbox{\rm any edge of a $2$-hammock is traversed by a through path.} 
\end{equation}

If $G$ is connected and $D \subseteq V(G)$ satisfies $|D|=k$ and $G-D$ is disconnected, then $D$ is called a $k$-\emph{disconnector}. The boundary of a $k$-hammock $H$ is a $k$-disconnector, unless $H$ is trivial.\\

\noindent
\scaps{Planar Hammocks.} The argument of Ma and Yu~\cite{MaYu1} for \sref{MaYuSep} asserts the following.

\begin{statement}\label{disc'}
Let $G$ be a $5$-connected nonplanar graph such that $G =G_{_1} \cup G_{_2}$, $|V(G_{_1} \cap G_{_2})| = 5$. If $|V(G)| > |V(G_{_2})| \geq 7$ and $G_{_2}$ has a an embedding in the closed disc with $V(G_{_1} \cap G_{_2})$ appearing on the boundary of the disc, then $TK_{_5} \subseteq G$.
\end{statement} 

\noindent
Indeed, the argument of Ma and Yu for \sref{MaYuSep} shows that the existence of a separation as in \sref{disc'} can be translated into a $TK_{_5}$, where the main engine is \cite[Theorem 4.3]{MaYu1}. 

A $5$-hammock $H$ of $G$ is called a \emph{planar hammock} if $|V(G)| > |V(H)| \geq 7$ and $H$ has an embedding in the closed disc with $\bnd{H}{}$ appearing on the boundary of the disc. 
For us it will be more convenient to use the following. 

\begin{statement}\label{disc}
A $5$-connected nonplanar graph containing a planar hammock contains a $TK_{_5}$. 
\end{statement}  

\noindent
\scaps{Blocks.} The \emph{blocks} of $G$ are its maximal $2$-connected components.
We refer to a single edge as being $2$-connected. 
The boundary of a block of $G$ is the set of cut vertices of $G$ contained in the block.

If $G$ is connected, then the blocks and cut vertices of $G$ define a unique tree called the \emph{block tree} of $G$; the leaves of which are blocks.  
A graph is called a \emph{chain} of blocks if its block tree is a (simple) path; and
is said to be a \emph{claw-chain} of blocks if its block tree has precisely $3$ leaves. By $xy$-\emph{chain} we mean a chain of blocks with each member of $\{x,y\}$ met by a leaf block of the chain. A subgraph of a chain that is also a chain is called a \emph{subchain}. 


\section{Proof of \sref{2con-main}}\label{proof-2con-main}
Suppose $G$ is a graph such that
\begin{equation}\label{K}
\mbox{$K_{_{2,3}} \iso K \subseteq G$, $V(K) = \{x_{_1},x_{_2},x_{_3},x_{_4},x_{_5}\}$
such that $d_{K}(x_{_1})=d_{K}(x_{_2})=3$.}
\end{equation}
We observe that if there is an $x_{_1}x_{_2}$-path $P$ in $G-\{x_{_3},x_{_4},x_{_5}\}$ such that 
\begin{enumerate}
\item[(F.1)] $G-P$ contains an $\{x_{_3},x_{_4},x_{_5}\}$-circuit, then $TK_{_5} \subseteq G$.
\item[(F.2)] Also, if there is a $v \in int P$ such that $G-P+v$ contains a $Z \cup \{v\}$-circuit $C$ satisfying $V(C) \cap \{x_{_3},x_{_4},x_{_5}\}=Z$ and $|Z|=2$, then $TK_{_5} \subseteq G$.
\end{enumerate}
We restate \sref{2con-main} as follows.

\begin{statement}\label{2con-1}
Let $G$ be a $5$-connected graph with $K\subseteq G$ as in (\ref{K}).
If $G$ has an induced $x_{_1} x_{_2}$-path $P$ satisfying $x_{_3},x_{_4},x_{_5} \notin V(P)$ and $\kappa(G-P) \geq 2$, then $TK_{_5} \subseteq G$.
\end{statement}

If $TK_{_5} \not\subseteq G$ and $G,P$ are as in \sref{2con-1}, then (F.1) fails in the sense that $G-P$ contains no $\{x_{_3},x_{_4},x_{_5}\}$-circuit. By \sref{frame-separation}, which is the characterization of $2$-connected graphs containing no circuit through some three prescribed elements (vertices or edges), the lack of such a circuit in $G-P$ is certified by what we refer to as an $\{x_{_3},x_{_4},x_{_5}\}$-frame (see \scaps{Definition A} below). 

\subsection{Frames}\label{frames}
A pair $(X,F)$ with $X\subset V(G)$ and $F\subseteq E(G-X)$ is called a \emph{separator} 
if the number of components of $G-X-F$ is greater than the \emph{capacity}

\begin{equation}
\gamma(X,F)\colon=|X|+\sum\left\{\lfloor\half|\bnd{Q}{G-X}|\rfloor\colon\;
                                    \mbox{\rm $Q$ is a component of $G[F]$}\right\}.
\end{equation}
Here, $G[F]=(V(F),F)$, where $V(F)$ is the set of vertices incident with $F$. Given a set $A\subseteq V(G) \cup E(G)$, we specify $(X,F)$ as an $A$-\emph{separator} if $A$ meets $> \gamma(X,F)$ components of $G-X-F$. A pair $(G,A)$ cannot admit both an $A$-circuit and an $A$-separator. In general, neither of these exists, except for quite special cases when 
the ``circuit-separator alternative'' holds. The simplest such case is  
the generalization of the well-known theorem of Whitney (case $A \subseteq V(G)$, $|A|=2$) proved by  Mesner and Watkins~\cite{mw}.

Throughout the remainder of this section $\kappa(G) =2$ and $U \subset E(G) \cup V(G)$, $|U|=3$. The $2$-hammocks of $G$ satisfying $|U \cap H|\le 1$ are referred to as
$U$-\emph{hammocks}. For two collections $\L$, $\L'$ of $U$-hammocks, let 
$\L \preceq \L'$ mean that each $H\in \L$ is a subgraph of some $H'\in \L'$, and let $\L$ 
be called {\sl maximal} if $\L \preceq \L'$ implies $\L=\L'$. In a maximal collection, a common vertex of two hammocks is an end of each. \\

\noindent
{\sc Definition A.}\\
A pair $\X=(X,F)$, where $X\subset V(G)$ and $F\subset E(G-X)$, is called a $U$-\emph{frame} if the subgraph $G(\X)\colon=(X\cup V(F),F)$ has exactly two components, and $G$ is the union of $G(\X)$ and a collection $\B$ of $U$-hammocks such that:\\
(A.1) a member of $\B$ has its ends in distinct components of $G(\X)$;\\
(A.2) $U \cap H \not= \emptyset$ holds for three $H \in \B$; and $u \in int H$ of some $H \in \B$ for each vertex $u \in U$.\\ 
(A.3) $0 \leq |X| \leq 2$; and  if $X=\emptyset$ then $|\B|=3$.\\
(A.4) $\B$ is maximal subject to (A.1-3).\\
The members of $\B$ will be referred to as $\X$-{\em hammocks}. Due to the maximality of $\B$, in case $F \not= \emptyset$:\\
(A.5) the ends of any two $\X$-hammocks in the same component of $G[F]$ are distinct;\\
(A.6) the components of $G[F]$ are $2$-connected and of order $\geq 3$.\\

\begin{figure}[htbp]
	\centering
		\scalebox{0.7}{\includegraphics{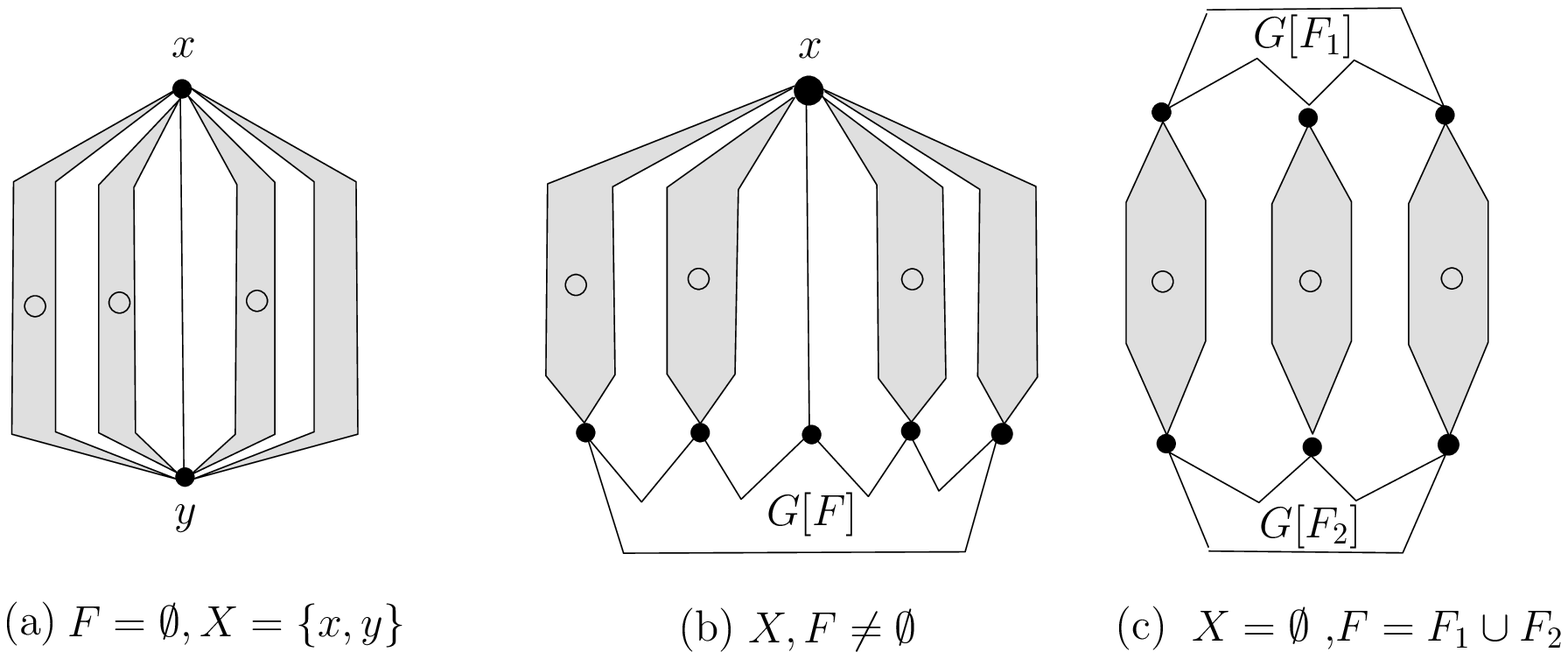}}
	\caption{Frames $\X=(X,F)$ separating $\{\circ,\circ,\circ\}$.}
	\label{frame}
\end{figure}

We say that a $U$-frame $\X$ \emph{separates} $U$.
Figure~\ref{frame} illustrates the $3$ typical forms of a frame.
The $\X$-hammocks meeting $U$ are called \emph{seeded} and the members of $U$ are called \emph{seeds}. A through path containing the seed of a seeded $\X$-hammock is called a \emph{seeded-through path}; the latter exists due to (\ref{through-path}). If $x$ is a vertex in a seeded $\X$-hammock, then a seeded-through path containing $x$ is called a \emph{seeded-through-$x$-path} (such paths do not necessarily exist). 

Since any nontrivial $\X$-hammock $H$ becomes $2$-connected once an ear connecting its ends is appended, $H$, if not $2$-connected, is a chain of blocks. Consequently, if $H$ is seeded with seed $x$, has ends $z,z'$, and $y \in int H -x$, then $H$ contains a $y\ell$-path $P$ satisfying $x \in int P$, $\ell \in \{z,z'\}$. Such a path is called a \emph{seeded-$y$-path}.   
The characterization of the $2$-connected graphs containing no circuit traversing $3$ prescribed elements of the graph (edges or vertices) is given by the following.   

\begin{statement}\label{frame-separation}
$G$ has no $U$-circuit if and only if $G$ has a $U$-frame. 
\end{statement}

Theorem \sref{frame-separation} is derivable from the (elegant) fact that: 
{\sl a $3$-connected (not necessarily simple) graph $(V,E)$ has no $A$-circuit, where $A\subset V\cup E$ and $|A|=3$, if and only if $A$ consists of edges, and $A$ is either an edge-cut, a claw ($K_{_{1,3}}$), or contains two parallel edges.} 

Ma and Yu~\cite{MaYu1} use a close variant of the notion of frames; their variant, they explain, can be derived from the Mesner-Watkins argument~\cite{mw}. 
The difference between Ma and Yu's variant to ours is that they ignore unseeded $U$-hammocks and essentially work with pairs $\X=(X,F)$ where $G(X)$ may have a single component instead of precisely $2$; also in their version, $U$ is always a set of vertices.  

\subsection{Proof of \sref{2con-1}} The reader should be reminded of the \scaps{agreement} specified in \S\ref{pre}. 
Throughout this section, $G$ is $5$-connected and contains no $TK_{_5}$, $P$ is as in \sref{2con-1}, $U = \{x_{_3},x_{_4},x_{_5}\}$ as these are defined in (\ref{K}), $v \in int P$ and for such a $v$, $U_{_v}$ denotes a set of size $3$ containing $v$ satisfying $|U \cap U_{_v}|=2$. Also, if $\X$ is a $Y$-frame in graph $J$ separating some $3$-element set $Y \subset V(J) \cup E(J)$, then the components of $J(\X)$ are referred to as $F_{_1}$ and $F_{_2}$, and the end of an $\X$-hammock $H$ contained in $F_{_i}$ is denoted $z_{_{i,H}}$, for $i=1,2$. 

To reach a contradiction, we show that a $U_{_v}$-circuit satisfying (F.2) exists. To this end, suffices to show that there is an 
\begin{equation}\label{enough}
\mbox{$\ell \ell'$-path $Q$ in $G-P$ satisfying $\ell,\ell' \in N_{G-P}(v)$ and $|V(Q) \cap U|=2$.}
\end{equation}  
Consequently, throughout this section, we assume, towards contradiction, that a path satisfying (\ref{enough}) does not exist. 

By (F.1) and \sref{frame-separation}, there is a $U$-frame in $G-P$. In what follows, we study $U$-frames in $G-P$ (under the assumption that a path satisfying (\ref{enough}) does not exist).\\

\noindent
\scaps{Reduction Steps.} Let $J$ be a graph with $\kappa(J) = 2$ and let $\X$ be a $W$-frame  for some $W \subset V(J), |W|= 3$. Let $H$ be a seeded $\X$-hammock with seed $x$. If $y \in V(H)$, then an $(X_{_y}=\{x,y,z_{_{1,H}} z_{_{2,H}}\})$-circuit in $R^1_{_H} = H+z_{_{1,H}} z_{_{2,H}}$, implies the existence of a seeded-through-$y$-path in $H$. Next, if $y,y' \in V(H)$, then an $(X_{_{y,y'}}=\{x,yy', z_{_{1,H}} z_{_{2,H}}\})$-circuit in $R^2_{_H} = R^2_{_H}(y,y') = H+z_{_{1,H}} z_{_{2,H}}+yy'$, implies the 
existence of a $yy'$-path in $G-P$ meeting precisely two seeds. Consequently, by \sref{frame-separation}:\\

\noindent
(R.1) \emph{If $y\in V(H)$ and $R^1_{_H}$ contains no seeded-through-$y$-path, then
there is an $X_{_y}$-frame in $R^1_{_H}$; and thus $y \in (int_{G-P} H) \sm \{x\}$.}\\

\noindent
(R.2) \emph{If $y,y' \in V(H)$ and $G-P$ does not contain a $yy'$-path meeting precisely two seeds, then there is an $X_{_{y,y'}}$-frame in $R^2_{_H}(y,y')$.}\\

For $i=1,2$, we refer to the frames defined in (R.i) as $\X_{_H} =(X_{_H},F_{_H})$.  The seeded $\X_{_H}$-hammocks containing $x,y,z_{_{1,H}} z_{_{2,H}}$ in (R.1), and $x,y y',z_{_{1,H}} z_{_{2,H}}$ in (R.2), are denoted by $K_{_H},K'_{_H},K''_{_H}$, respectively. By $\P_{_H}$ we mean a $2$-linkage in $G-P$ connecting $\{z_{_{1,H}}, z_{_{2,H}}\}$ to the ends of $K''_{_H}$; such a linkage is clearly contained in $K''_{_H} - z_{_{1,H}} z_{_{2,H}}$. Let $\ell_{_{1,H}},\ell_{_{2,H}}$ be the ends of $K''_{_H}$ linked to $z_{_{1,H}}, z_{_{2,H}}$ in $\P_{_H}$, respectively. Call the components of $R^i_{_H}(\X_{_H})$, $i=1,2$, containing $\ell_{_{1,H}},\ell_{_{2,H}}$, by $F_{_{1,H}}$ and $F_{_{2,H}}$, respectively.\\

\noindent
\scaps{Constructions.} We list several settings and constructions of paths therein which yield a path satisfying (\ref{enough}) and that rise frequently in subsequent arguments.
Let $J$ be a graph with $\kappa(J) = 2$ and let $\X$ be a $W$-frame  for some $W \subset V(J), |W|= 3$.\\

\noindent
(S.1) \emph{Let $H,H'$ be seeded $\X$-hammocks and let $y \in int_J H$.   
Let $\X_{_{H'}}$ be an $X_{_u}$-frame 
or an $X_{_{u,u'}}$-frame in $R^i_{_{H'}}$ per (R.i) for $i=1,2$, where $u,u' \in V(H')$.
If $H$ contains a seeded-through-$y$-path, then $J$ contains a $(y,\{u,u'\})$-path meeting precisely two seeds.} \\

\proof
Let $S$ be a seeded-through-$y$-path in $H$ such that $z_{_{i,H}},y,x,z_{_{3-i,H}}$ appear in this order along $S$ for some $i=1,2$, and where $x$ is the seed of $H$. A path satisfying the assertion is contained in the union of the following paths: the $y z_{_{3-i,H}}$-subpath of $S$, a path in $F_{_{3-i}}$ linking the ends of $H$ and $H'$, the the $\ell_{_{3-i,H'}} z_{_{3-i,H'}}$-path in $\P_{_{H'}}$, a path in $F_{_{3-i,H'}}$ linking the ends of $K_{_{H'}}$ and $K''_{_{H'}}$, a seeded-through path in  $K_{_{H'}}$, a path in $F_{_{i,H'}}$ linking the ends of $K_{_{H'}}$ and $K'_{_{H'}}$, a path in $K'_{_H} -uu'$ linking $z_{_{i,K'_{_{H'}}}}$ to a member of $\{u,u'\}$ not meeting the other end of $K'_{_H}$. Easy to see that the last path exists.\QED \\

\noindent
(S.2) Let $H,H'$ be seeded $\X$-hammocks and let $y \in V(H), y' \in V(H')$.
Let $\X_{_H}$ and $\X_{_{H'}}$ be an $X_{_y}$-frame and an $X_{_{y'}}$-frame in $R^1_{_H}$ and $R^1_{_{H'}}$, respectively. Pick an $i \in \{1,2\}$ and consider the following union of paths:
a $y z_{_{i,K_{_H}}}$-path in $K'_{_H}$ not meeting $z_{_{3-i,K_{_H}}}$ (such a path exists by (\ref{through-path}) and $y \in int_J H$ by (R.1)), a path in $F_{_{i,H}}$ linking the ends of $K'_{_H}$ and $K''_{_H}$, the $\ell_{_{i,H}} z_{_{i,H}}$-path in $\P_{_H}$, a path in $F_{_i}$ linking the ends of $H$ and $H''$, where $H''$ is the third seeded $\X$-hammock, a seeded-through path in $H''$, a path in $F_{_{3-i}}$ linking the ends of $H'$ and $H''$, the $z_{_{3-i,H'}} \ell_{_{3-i,H'}}$-path in $\P_{_{H'}}$, a path in $F_{_{3-i,H'}}$ connecting the ends of $K_{_{H'}}$ and $K''_{_{H'}}$, a seeded-through path in $K_{_{H'}}$, a path in $F_{_{i,H'}}$ linking the ends of $K_{_{H'}}$ and $K'_{_{H'}}$, a $z_{_{i,K'_{_{H'}}}} y'$-path in $K'_{_{H'}}$ not meeting
$z_{_{i,K'_{_{H'}}}}$. 

The resulting subgraph $R$ is called a \emph{trail}. Such is a $yy'$-path in $J$ meeting precisely two seeds provided that 
$R \cap F_{_i}$ and $R \cap F_{_{i,H'}}$ are disjoint. By the definition of frames, this does not occur provided $F_{_i}$ and $F_{_{i,H'}}$ coincide into a singleton.  

Due to symmetry, the above construction defines $4$ distinct trails: indeed we may choose to start the construction from either $y$ or $y'$ (here we started from $y$) for each such choice we may choose $i \in \{1,2\}$.\\ 

\noindent
(S.3) \emph{Let $y \in int_J H$, $y' \in int_J H'$, where $H,H'$ are a seeded and an unseeded $\X$-hammocks, respectively. If $H$ contains a seeded-through-$y$-path, then $J$ contains a $yy'$-path meeting precisely two seeds.}\\

\proof
Let $S$ be a seeded-through-$y$-path in $H$ such that $z_{_{i,H}},y,x,z_{_{3-i,H}}$ appear in this order along $S$ for some $i=1,2$, and where $x$ is the seed of $H$. A path satisfying the assertion is contained in the union of the following paths: the $y z_{_{3-i,H}}$-subpath of $S$, a path in $F_{_{3-i}}$ linking the ends of $H$ and $H''$, where $H''$ is some seeded $\X$-hammock, a seeded-through path in $H''$, a path in $F_{_i}$ linking the ends of 
$H''$ and $H'$, a $z_{_{i,H'}} y'$-path in $H'$ not meeting $z_{_{3-i,H'}}$.\QED \\

\noindent
(S.4) Let $y \in V(H)$, $y' \in int_J H'$, where $H$ and $H'$ are a seeded and unseeded $\X$-hammocks, respectively. Let $\X_{_H}$ be an $X_{_y}$-frame in $R^1_{_H}$.
Pick an $i \in \{1,2\}$ and consider the following union of paths:
a $y z_{_{i,K'_{_H}}}$-path in $K'_{_H}$ not meeting $z_{_{3-i,K'_{_H}}}$, a path in $F_{_{i,H}}$ linking the ends of $K_{_H}$ and $K'_{_H}$,
a seeded-through path in $K_{_H}$, a path in $F_{_{3-i,H}}$ linking the ends of $K_{_H}$ and $K''_{_H}$, the $\ell_{_{3-i,H}} z_{_{3-i,H}}$-path in $\P_{_H}$, a path in $F_{_{3-i}}$ linking the ends of $H$ and $H''$ where $H''$ is some seeded $\X$-hammock, a seeded-through path in $H''$, a path in $F_{_i}$ linking the ends of $H''$ and $H'$, a $z_{_{i,H'}} y'$-path in $H'$ not meeting $z_{_{3-i,H'}}$. 

The resulting subgraph $R$ is called a \emph{poor trail}. Such is a $yy'$-path in $J$ meeting precisely two seeds provided that 
$R \cap F_{_i}$ and $R \cap F_{_{i,H}}$ are disjoint. By the definition of frames, this does not occur provided $F_{_i}$ and $F_{_{i,H}}$ coincide into a singleton. 

Due to symmetry, the above construction defines $2$ distinct poor trails as we may choose $i \in \{1,2\}$.\\

\noindent
\scaps{Properties of $U$-frames in $G-P$.} Let $\X$ be a $U$-frame in $G-P$, $\B$ its collection of $\X$-hammocks, and let $v \in int  P$. We consider several properties of $\X$ and $N_{G-P}(v)$. \\

\noindent
(\sref{2con-1}.A) \emph{$N_{G-P}(v) \subseteq \bigcup_{H \in \B} V(H)$.}\\

\proof 
The following two claims (\sref{2con-1}.A.1) and (\sref{2con-1}.A.2) are sufficient.\\

\noindent
(\sref{2con-1}.A.1) \emph{If $|N_G(v) \cap V(F_{_i})| \geq 2$, for some $i=1,2$, then a $U_{_v}$-circuit satisfying (F.2) exists.}\\

\proof
Let $y,y' \in N_G(v) \cap V(F_{_i})$, for some $i=1,2$. Let $\P$ be a $2$-linkage
in $F_{_i}$ connecting $\{y,y'\}$ and two ends in $F_{_i}$ of some two seeded $\X$-hammocks, namely $H$ and $H'$. Let $Y$ and $Y'$ be seeded-through paths in $H$ and $H'$, respectively. Let $Q$ be a path in $F_{_{3-i}}$ connecting the ends of $H$ and $H'$ in $F_{_{3-i}}$. Clearly, $\P \cup Y\cup Y' \cup Q \cup \{yv,y'v\}$ is a $U_{_v}$-circuit satisfying (F.2). \inQED \\

By (\sref{2con-1}.A.1), and since $|N_{G-P}(v)| \geq 3$ (as $P$ is induced), and since $G(\X)$ has precisely two components, $N_H(v) \not= \emptyset$ for at least one $H \in \B$. To conclude then, we prove the following. \\

\noindent
(\sref{2con-1}.A.2) \emph{If $N_G(v) \cap int_{G-P} F_{_i} \not= \emptyset$, for some $i=1,2$, and $N_H(v) \not=\emptyset$, for some $H \in \B$, then 
a $U_{_v}$-circuit satisfying (F.2) exists.}\\

\proof 
Here we prefer to use the following notation. 
Let $H_{_i}$ denote the seeded $\X$-hammocks and let $z^j_{_i}$ denote the end of $H_{_i}$ in $F_{_j}$, for $i=1,2,3$ and $j=1,2$. Let $S_{_i}$ denote a seeded-through path in $H_{_i}$, and let $Q_{_{i,x}}$ denote a seeded-$x$-path in $H_{_i}$, where $x \in N_{H_{_i}}(v)$. 
By $P^i_{_{k,\ell}}$ we mean a $z^i_{_k} z^i_{_{\ell}}$-path in $F_{_i}$, for $i=1,2$ and $k,\ell=1,2,3$, $k \not= \ell$. Let $y \in N_{G}(v) \cap int_{G-P} F_{_r}$, for some $r=1,2$. By $F_{_{k,\ell}}$ we mean a $(y, \{z^r_{_k},z^r_{_{\ell}}\})$-$2$-fan in $F_{_r}$, for $r=1,2$, $k,\ell=1,2,3$, $k \not= \ell$. Let $\P_{_{j,k,\ell}}$ denote a $(\{y,z^r_{_j}\},\{z^r_{_k},z^r_{_{\ell}}\})$-$2$-linkage in $F_{_r}$. We observe the following.  

\noindent
\bfm{(a)} If $x \in N_{H_{_j}}(v)$, $j=1,2,3$, then $Q_{_{j,x}}$ ends in $z^r_{_j}$. 
Otherwise, a path satisfying (\ref{enough}) is contained in $Q_{_{j,x}} \cup F_{_{j,k}} \cup
S_{_k} \cup P^{3-r}_{_{j,k}}$, where $j \not= k$, $j,k =1,2,3$.

\noindent
\bfm{(b)} By \bfm{(a)} and (\sref{2con-1}.A.1), if $N_{H_{_j}}(v) \not= \emptyset$, for some $j=1,2,3$, then $N_{F_{_{3-r}}}(v) \cap int_{G-P} F_{_{3-r}} = \emptyset$. 

\noindent
\bfm{(c)} If $x \in N_{H_{_j}}(v)$, then $N_{H_{_k}}(v) = N_{H_{_{\ell}}}(v) = \emptyset$, where $\{j,k,\ell\} = \{1,2,3\}$. To see this, let $z \in N_{H_{_k}}(v)$.
By \bfm{(a)}, $Q_{_{j,x}}$ and $Q_{_{k,z}}$ end at $z^r_{_j}$ and $z^r_{_k}$, respectively.
Thus, a path satisfying (\ref{enough}) is contained in $Q_{_{j,x}} \cup Q_{_{j,z}} \cup P^r_{_{j,k}}$.

\noindent
\bfm{(d)} If $x,z \in N_{H_{_j}}(v)$, $j=1,2,3$, then a path satisfying (\ref{enough}) exists.
Let $\P \subseteq H_{_j}$ be an $(\{x,z\},\{z^r_{_j},z^{3-r}_{_j}\})$-$2$-linkage.
The seed of $H_{_j}$ is not contained in $\P$. Indeed, if so, a path satisfying (\ref{enough}) is contained in $\P \cup P^r_{_{j,k}} \cup P^{3-r}_{_{j,k}} \cup S_{_k}$, where $j\not= k =1,2,3$. Consequently, a path satisfying (\ref{enough}) is contained in $R \cup \P_{_{j,k,\ell}} \cup S_{_k} \cup P^{3-r}_{_{k,\ell}} \cup S_{_{\ell}}$, where $R \in \P$ and has $z^r_{_j}$ as an end and $\{j,k,l\}=\{1,2,3\}$.

\noindent
\bfm{(e)} If $x \in N_{G-P}(v) \cap int_{G-P} H$, where $H\in \B$ and unseeded, then a path satisfying (\ref{enough}) exists. To see this, let $R\subset H$ be an $xw$-path such that 
$w$ is the end of $H$ in $F_{_r}$ and such that $R$ meets only one end of $H$; by (\ref{through-path}), $R$ exists. Let $\P$ be a $(\{w,y\},\{z^r_{_k},z^r_{_{\ell}}\})$-$2$-linkage in $F_{_r}$. A path satisfying (\ref{enough}) is contained in $R \cup \P \cup S_{_\ell} \cup S_{_k} \cup P^{3-r}_{_{k,\ell}}$, where $k,\ell=1,2,3$, $k\not=\ell$.\\

To conclude the proof consider the following. We have that $N_H(v) \not= \emptyset$ for some $H \in \B$, by (\sref{2con-1}.A.1). Presence of an unseeded $\X$-hammock implies that at least one of the components of $G'(\X)$, where $G' = G-P$, is a singleton. In our setting, $F_{_{3-r}}$ will be a singleton. 
Thus, by (\sref{2con-1}.A.1) and (b), (e), we may assume
that $H$ is seeded, regardless of the existence of unseeded $\X$-hammocks. Now, as $N_{F_{_{3-r}}(v)} \cap int_{G-P} F_{_{3-r}}= \emptyset$ and $|N_{F_{_r}}(v)| =1$, by (b) and (\sref{2con-1}.A.1), respectively, then either $|N_H(v)| \geq 2$ contradicting (d), or $N_{G-P}(v)$ meets at least two seeded hammocks contradicting (c), or $N_{G-P}(v)$ meets the interior of an unseeded $\X$-hammock contradicting (e).\inQED \\
\QED \\

\noindent
(\sref{2con-1}.B) \emph{If $y \in N_{G-P}(v) \cap int_{G-P}H$ where $H$ is a seeded $\X$-hammock containing a seeded-through-$y$-path,
then no other seeded $\X$-hammock meets $\geq 2$ members of $N_{G-P}(v)$.}\\

\proof
Let $u,u'\in N_{H'}(v)$, where $H' \not= H$ is a seeded $\X$-hammock. Clearly $u,u' \not= y$.  
By (R.2), there is an $X_{_{u,u'}}$-frame in $R^2_{_{H'}}(u,u')$. Thus a path satisfying 
(\ref{enough}) exists, by (S.1). \QED \\

\noindent
(\sref{2con-1}.C) \emph{Let $y \in N_{H}(v)$, $y' \in N_{G-P}(v) \cap int_{G-P} H'$, where $H,H'$ are seeded $\X$-hammocks. If there is no seeded-through-$y$-path in $H$ and $H'$ has a seeded-through-$y'$-path, then a path satisfying (\ref{enough}) exists.}\\

\proof
By (R.1), there is an $X_{_y}$-frame in $R^1_{_{H}}$. Thus a path satisfying (\ref{enough}) exists, by (S.1). \QED \\

\noindent
\scaps{Properties of minimum $U$-frames.} For a $U$-frame $\X$ in $G-P$, put
\begin{equation}\label{para1}
\alpha(\X) = |\bigcup \{V(H): \mbox{$H$ is a seeded $\X$-hammock}\}|.
\end{equation}
We write $\X_{_1} \preceq \X_{_2}$ if $\alpha(\X_{_1}) \leq \alpha(\X_{_2})$, where $\X_{_1},\X_{_2}$ are $U$-frames in $G-P$. We say that a $U$-frame $\X$ is {\sl minimum}
if $\X' \preceq \X$ implies $\X' = \X$ for every $U$-frame $\X'$. \\

\noindent
(\sref{2con-1}.D) \emph{Let $\X$ be a minimum $U$-frame. If $y \in N_{H}(v)$, $y' \in N_{H'}(v)$, where $H,H'$ are seeded $\X$-hammocks, then $H,H'$ contain a seeded-through-$y$-path and a seeded-through-$y'$-path, respectively.}\\

\proof
At least one of $y,y'$ is contained internally in $H,H'$, respectively; otherwise the claim follows trivially from (\ref{through-path}). Consequently, by (\sref{2con-1}.C), we may assume that $H,H'$ contain no seeded-through-$y$-path and no seeded-through-$y'$-path, respectively. Thus, by (R.1), there are an $X_{_y}$-frame and an $X_{_{y'}}$-frame in $R^1_{_H}$ and $R^1_{_{H'}}$, respectively.  

The assumption that there is no path satisfying (\ref{enough}) implies that all (symmetric) trails produced by (S.2) fail to provide such a path. This occurs provided  $F_{_i},F_{_{i,H}}, F_{_{i,H'}}$ all coincide into a singleton $q_{_i}$, for $i=1,2$. 
The pair $\X'=(\{q_{_1},q_{_2}\}, \emptyset)$ then is a $U$-frame with a collection $\B'$ of
$U$-hammocks satisfying $K_{_H},K_{_{H'}},H'' \in \B'$ implying that $\alpha(\X') < \alpha(\X)$ since $K_{_H} \subset H,K_{_{H'}} \subset H'$. \QED \\
 
\noindent
(\sref{2con-1}.E) \emph{Let $\X$ be a minimum $U$-frame. If $N_H(v),N_{H'}(v) \not= \emptyset$ and $N_H(v) \not= N_{H'}(v)$ for some two seeded $\X$-hammocks $H,H' \in \B$, then $N_{G-P}(v) \subseteq V(H) \cup V(H')$.}\\

\proof
Let $y \in N_H(v), y' \in N_{H'}(v), y\not= y'$. Let $x,x'$ be the seeds of $H,H'$, respectively. Seeded-through-$b$-paths, $b \in \{y,y'\}$ in $H$ and $H'$, respectively, exist by (\sref{2con-1}.D). Such have noncorresponding orderings in the sense that:
if $H$ has a seeded-through-$y$-path such that $z_{_{i,H}},y,x,z_{_{3-i,H}}$ appear in this order along the path, then a seeded-through-$y'$-path in $H'$ has $z_{_{i,H'}},x',y',z_{_{3-i,H'}}$ appear in this order along it. Indeed, if not so, and the orderings correspond, then since $y\not= y'$, a path satisfying (\ref{enough}) clearly exists. Consequently, $N_{G-P}(v)$ does not meet the interior of the third seeded $\X$-hammock. 

Let then $H'' \in \B$ be unseeded such that $y'' \in N_{G-P}(v) \cap int_{G-P}H''$. 
By the definition of frames, $F_{_i}$ is a singleton for some $i=1,2$. Consequently, at least one of $\{y,y'\}$, say $y$, is contained internally in its corresponding $\X$-hammock. A path satisfying (\ref{enough}) then exists by (S.3). Such an $H''$ exists by (\sref{2con-1}.A) and $F_{_i}$ being a singleton for some $i=1,2$.\QED \\

\noindent
(\sref{2con-1}.F) \emph{Let $\X$ be a minimum $U$-frame. If $N_{G-P}(v) \cap int_{G-P}H \not= \emptyset$ for some seeded $\X$-hammock $H$, then $N_{G-P}(v) \subseteq V(H)$.}\\

\proof
By (\sref{2con-1}.E), if $N_{G-P}(v) \cap int_{G-P}H' \not= \emptyset$ for some seeded $\X$-hammock distinct of $H$, then $N_{G-P}(v) \subseteq V(H) \cup V(H')$. This implies that  $H,H'$ meet at least three distinct members of $N_{G-P}(v)$. A path satisfying (\ref{enough}) then exists, by (\sref{2con-1}.B) and (\sref{2con-1}.D).

Let $H'$ be unseeded then. Thus, $H$ has no seeded-through-$y$-path, $y \in N_{G-P}(v) \cap int_{G-P}H$; otherwise a path satisfying (\ref{enough}) exists by (S.3). Thus, by (R.1), there is an $X_{_y}$-frame in $R^1_{_H}$. The assumption that no path satisfying (\ref{enough}) implies that all (symmetric) poor trails produced by (S.4) fail to provide such a path. 
This occurs provided $F_{_i}$ and $F_{_{i,H}}$ coincide into a singleton $q_{_i}$ for $i=1,2$. Consequently, $\X'=(\{q_{_1},q_{_2}\}, \emptyset)$ is a $U$-frame with a collection $\B'$ of $U$-hammocks satisfying $K_{_H},H'',H''' \in \B'$, where $H'',H'''$ are the remaining seeded $\X$-hammocks. As $K_{_H} \subset H$, it follows that $\alpha(\X') < \alpha(\X)$\QED \\ 

\noindent
\scaps{Consequences.} We learn the following from minimum $U$-frames in $G-P$.\\
 
\noindent 
(\sref{2con-1}.G) \emph{A seed is not adjacent to $int P$.}\\

\proof
Let $x \in U$ be adjacent to $int P$ and let $\X$ be a minimum $U$-frame with $H$ the $\X$-hammock containing $x$. Let $v \in int P$ such that $x \in N_{H}(v)$. As $x \in  int_{G-P} H$, by definition of frames, we have that $N_{G-P}(v) \subseteq V(H)$, by (\sref{2con-1}.F). Let $y \in N_{H}(v) \sm \{x\}$. A $\{yx, z_{_{1,H}} z_{_{2,H}} \}$-circuit clearly exists in 
$H+\{yx, z_{_{1,H}} z_{_{2,H}} \}$. Such a circuit is easily extended into a $yx$-path satisfying (\ref{enough}).\QED \\

Property (\sref{2con-1}.G) implies that the seeded $\X$-hammocks of a $U$-frame $\X$ have order $\geq 4$ each. Thus, $5$-connectivity then asserts that the interior of such hammocks must be adjacent to $int P$. Property (\sref{2con-1}.F) then implies that if $H_{_i}$, $i=1,2,3$, are the seeded $\X$-hammocks of a minimum $\X$, then there is a partition $L_{_1} \cup L_{_2} \cup L_{_3} \cup L_{_4} = int P$ such that $L_{_i}$, $i=1,2,3$, consists of the vertices of $int P$ adjacent to $int_{G-P}H_{_i}$, and $L_{_4}$ is the complement of $L_{_1} \cup L_{_2} \cup L_{_3}$. We have seen that for $i=1,2,3$, $L_{_i} \not= \emptyset$; moreover, if $v \in L_{_i}$, then $N_{G-P}(v) \subseteq V(H_{_i})$, by (\sref{2con-1}.F).

Let $v \in L \in \{L_{_1},L_{_2},L_{_3}\}$ and let $H$ denote the seeded $\X$-hammock corresponding to $L$. Clearly, (R.2) applies to $H$ and any two members of $N_H(v)$; or a path satisfying (\ref{enough}) exists. Let $\L_{_{v,H}}$ denote the $\X_{_H}$-frames of type (R.2) in $R^2_{_H}$ where (R.2) is applied to the seed of $H$ and two members of $N_H(v)$; we have just seen that $\L_{_{v,H}} \not= \emptyset$. For $y \in V(H)$, that is not the seed of $H$, let $\L_{_{y,H}}=\emptyset$ if $H$ contains a seeded-through-$y$-path. Otherwise let $\L_{_{y,H}}$ denote the $\X_{_H}$-frames of type (R.1) in $R^1_{_H}$, where (R.1) is applied to $y$ and the seed of $H$. Put $\L_{_H} = (\bigcup_{v \in L} \L_{_{v,H}}) \cup (\bigcup_{y\in V(H)} \L_{_{y,H}}$). 
For $\X_{_H} \in \L_{_H}$, put $\beta(\X_{_H}) = |V(K_{_H})|$ and let $\beta = \min \{\beta(\X_{_H}): \X_{_H} \in \L_{_H}\}$. \\ 

We are now ready to prove \sref{2con-1} and consequently \sref{2con-main}.

\noindent
\bfm{Proof of \sref{2con-1}.} Let $\X$ be a minimum $U$-frame in $G-P$ and let $L_{_i}$, $1 \leq i\leq 4$ be the partition of $int P$ induced by $\X$ as defined above.
Let $L \in \{L_{_1},L_{_2},L_{_3}\}$, let $H$ be the seeded $\X$-hammock corresponding to $L$, and let $x$ be its seed.   
\begin{equation}\label{choice'}
\mbox{Choose $\X_{_H} \in \L_{_H}$ satisfying $\beta(\X_{_H})=\beta$.}
\end{equation}
The minimality of $\X_{_H}$ as in (\ref{choice'}) implies that for every 
$v' \in L$ and $y \in V(K_{_H})\sm \{x\}$,
\begin{equation}\label{props}
\mbox{(i) $|N_{K_{_H}}(v')| \leq 1$ and (ii) $K_{_H}$ contains a seeded-through-$y$-path.}
\end{equation}
{\sl Proof.}
To see (\ref{props}), let $\ell, \ell'$ denote the ends of $K_{_H}$. Observe that an $\{x,yy',\ell \ell'\}$-circuit in $K_{_H} + \{yy', \ell \ell'\}$, for some two $y,y' \in N_{K_{_H}}(v')$, easily extends into a $yy'$-path satisfying (\ref{enough}). Consequently, if (\ref{props}) (i) is not satisfied, then an $\{x,yy',\ell \ell'\}$-frame $\X_{_{K_{_H}}}$  exists in $K_{_H} + \{yy', \ell \ell'\}$, by (R.2). In a similar manner, if $K_{_H}$ contains no seeded-through-$y$-path, then $K_{_H} +\ell \ell'$ has an $\{x,y,\ell \ell'\}$-frame $\X_{_{K_{_H}}}$ by in (R.1). 
  
Assume now, to the contrary, that (\ref{props}) (i) or (ii) is not satisfied for some $v' \in L$ and $u,u'\in N_{K_{_H}}(v')$ or some $y \in V(K_{_H}) \sm \{x\}$, respectively, and let $\X_{_{K_{_H}}}$ be a frame in $R^2_{K_{_H}}(u,u')$ and $R^1_{K_{_H}}$, respectively. Consider the following frame $\Y$. Let $G' \in \{R^1_{K_{_H}},R^1_{K_{_H}}(u,u')\}$ and put: $G'(\Y) = G'(\X_{_{K_{_H}}}) = F_{_{1,K_{_H}}} \cup F_{_{2,K_{_H}}}$, and the $\Y$-hammocks are $K_{_{K_{_H}}},K'_{_{K_{_H}}}$, and the third $\Y$-hammock is the union of $K''_{_{K_{_H}}}$, $F_{_{1,H}}, F_{_{2,H}}, K'_{_H}$, and $K''_{_H}$. 

If $\X_{_{K_{_H}}}$ is the result of a violation of (\ref{props}) (i) for some $v' \in L$ and $u,u' \in N_{K_{_H}}(v')$, then $\Y \in \L_{_{v',H}}$ (and thus in $\L_{_H}$).
If $\X_{_{K_{_H}}}$ is the result of a violation of (\ref{props}) (ii) for some $y \in V(K_{_H}) \sm \{x\}$, then $\Y \in \L_{_{y,H}}$ (and thus in $\L_{_H}$).

Property (\ref{props}) now follows since $K_{_{K_{_H}}} \subset K_{_H}$ implying that $\beta(\Y) < \beta(\X_{_H})$; contradicting (\ref{choice'}). \inQED \\

By (\sref{2con-1}.G), $|V(K_{_H})| \geq 4$; thus, $5$-connectivity and (\ref{props}) then imply that there are a $v'$ and a $y$ satisfying $v' \in L$ and a $\{y\} = int_{G-P} K_{_H} \cap N_{H}(v')$. As $N_{G-P}(v') = N_{H}(v')$, then $N_{H}(v') \sm \{y\} \subseteq H-int_{G-P} K_{_H}$, by (\sref{2con-1}.F), (\sref{2con-1}.G) and (\ref{props}). 

To conclude the proof, suffices that we show that 
\begin{equation}\label{enough-2}
\mbox{there are $i=1,2$ and $y',y''$ such that $y',y'' \in V(F_{_{i,H}}) \cap (N_H(v') \sm \{y\})$.} 
\end{equation}
Indeed, if so, it is easy to see that: a $(\{y',y''\},\{z_{_{i,K_{_H}}}, \ell_{_{i,H}}\})$-$2$-linkage in $F_{_{i,H}}$, the linkage $\P_{_H}$, a seeded-through path in $K_{_H}$, and a seeded-through path in $H'$, where $H'$ is some seeded $\X$-hammock, can be extended into a $y'y''$-path satisfying (\ref{enough}).

To establish (\ref{enough-2}), we argue as follows. Let $y' \in N_H(v') \sm \{y\}$. 
Let $S \subseteq K_{_H}$ be the seeded-through-$y$-path in $K_{_H}$, which exists by (\ref{props}), such that $z_{_{i,K_{_H}}},y,x,z_{_{3-i,K_{_H}}}$ appear in this order along $S$, where $x$ is the seed of $K_{_H}$. Suffices that we argue the following:\\
\bfm{(a)} $y' \notin F_{_{i,H}} \cup (K'_{_H} - z_{_{3-i,K'_{_H}}})$. Indeed, if this is not so, then a path satisfying (\ref{enough}) is contained in the union of the following paths: the $y z_{_{3-i,K_{_H}}}$-subpath of $S$, a path in $F_{_{3-i,H}}$ linking the ends of $K_{_H}$ and $K''_{_H}$, the members of $\P_{_H}$, a path in $F_{_{3-i}}$ linking the ends of $H$ and $H'$, where $H'$ is some seeded $\X$-hammock, a seeded-through path in $H'$, a path in $F_{_i}$ linking the ends of $H'$ and $H$, a path in $F_{_{i,H}} \cup (K'_{_H} - \{z_{_{3-i,K'_{_H}}}\})$ linking $y'$ with $\ell_{_{i,H}}$.

\noindent    
\bfm{(b)} $y' \notin int_{G-P} K''_{_H}$. If not so, then a $(\{y', \ell_{_{3-i,H}}\}, \{z_{_{1,H}},z_{_{2,H}}\})$-$2$-linkage is clearly contained in $K''_{_H} - z_{_{1,H}}z_{_{2,H}}$ and is easily extended into a path satisfying (\ref{enough}).\QED

\section{Apex graphs containing $\mbs{K_{_{2,3}}}$}\label{apex-main}
To outline our approach for proving \sref{k23}, suppose now that $G$ is a $5$-connected nonplanar apex graph containing no $TK_{_5}$ and $K \subseteq G$, where $K$ is as in (\ref{K}). By \sref{MaYuk4}, we may assume throughout this section that
\begin{equation}\label{stable}
\mbox{$\{x_{_1},x_{_2}\}$  and  $\{x_{_3},x_{_4},x_{_5}\}$ are stable.}
\end{equation}

By \sref{always}, there is an induced $x_{_1} x_{_2}$-path $P$ such that $x_{_3},x_{_4},x_{_5}\notin V(P)$ and with a $P$-bridge meeting at least two $2$-valent vertices of $K$.
We propose that two cases be examined:\\
\scaps{Connected Case:} $P$ can be chosen so that there is a $P$-bridge containing $\{x_{_3},x_{_4},x_{_5}\}$, or\\
\scaps{Disconnected Case:} such a choice is impossible.

In the \scaps{Connected Case}, by \sref{rooted-nonseparating-path} with $Y=\{x_{_3},x_{_4},x_{_5}\}$ and $X =\emptyset$ ($X,Y$ - see formulation of \sref{rooted-nonseparating-path}), $P$ can be chosen such that 
\begin{equation}\label{P}
\mbox{$G-P$ is connected with each leaf block meeting $\{x_{_3},x_{_4},x_{_5}\}$ internally.} 
\end{equation}
By \sref{2con-main} we may assume that $\kappa(G-P) = 1$ and thus $G-P$ is a chain or a claw-chain of blocks. By \sref{notapex}, assuming that $TK_{_5} \not\subseteq G$ implies that neither $x_{_1}$ nor $x_{_2}$ is an apex vertex of $G$. To reach contradiction in this case we that

\begin{statement}\label{1con-1}
if $\kappa(G-P) = 1$, each leaf block of $G-P$ meets $\{x_{_3},x_{_4},x_{_5}\}$ internally, and $a \in V(G)\sm\{x_{_1},x_{_2}\}$ is an apex vertex of $G$, then $TK_{_5} \subseteq G$. 
\end{statement}

For the \scaps{Disconnected Case} we argue the following. 
By \sref{always}, \sref{2con-main} and \sref{1con-1}, we may assume, without loss of generality, that $P$ can be chosen so that there is $P$-bridge containing $\{x_{_4},x_{_5}\}$, and there is no such $P$ with a $P$-bridge containing $\{x_{_3},x_{_4},x_{_5}\}$. Thus, by \sref{rooted-nonseparating-path} with $Y  =\{x_{_4},x_{_5}\}$ and $X = \{x_{_3}\}$, $P$ can be chosen so that $G-P$ consists of precisely two components: 
\begin{equation}\label{P-disconnected}
\mbox{$x_{_3}$ and a chain of blocks $H$ with each leaf block meeting $\{x_{_4},x_{_5}\}$ internally.}
\end{equation} 
If $\kappa(H) \geq 2$, we utilize the assumption that (F.2) fails and use \sref{frame-separation} to prove that

\begin{statement}\label{2con-2}
if $\kappa(H) \geq 2$ and $a \in V(G) \setminus \{x_{_1},x_{_2}\}$ is an apex vertex of $G$, then $TK_{_5} \subseteq G$. 
\end{statement}

\noindent
To conclude our arguments for the \scaps{Disconnected Case}, and indeed our proof of \sref{k23}, we show that 

\begin{statement}\label{1con-2}
if $\kappa(H) = 1$ and $a \in V(G) \setminus \{x_{_1},x_{_2}\}$ is an apex vertex of $G$, then $TK_{_5} \subseteq G$. 
\end{statement}

In the remainder of this section we prove \sref{1con-1}, \sref{2con-2}, and
\sref{1con-2}. 


\subsection{Rooted nonseparating paths}\label{rooted-paths} 
The purpose of this section is to prove \sref{rooted-nonseparating-path} and \sref{always}; usage of these is outlined above. 
Lemma \sref{rooted-nonseparating-path} is modeled after \cite[Lemma 2.3]{MaYu1} and serves as a generalized variant of \cite[Lemma 2.3]{MaYu1}. Our proof of \sref{rooted-nonseparating-path} follows from the three claims: (\sref{rooted-nonseparating-path}.A), (\sref{rooted-nonseparating-path}.B), and (\sref{rooted-nonseparating-path}.C) listed below. The proof of (\sref{rooted-nonseparating-path}.A) is precisely the proof of an analogous claim that Ma and Yu use to prove \cite[Lemma 2.3]{MaYu1}; we do not see a need to change their argument. To prove (\sref{rooted-nonseparating-path}.B), we use the approach of Ma and Yu for an analogous claim in their argument, but the implementation of their approach here requires adjustments resulting from the differences between \sref{rooted-nonseparating-path} and \cite[Lemma 2.3]{MaYu1}. Claim (\sref{rooted-nonseparating-path}.C) is new. 
Consequently, we attribute the following to Ma and Yu, but see no option but to include a complete proof of \sref{rooted-nonseparating-path}.

\begin{statement}\label{rooted-nonseparating-path}
Suppose $\kappa(G) \geq 5$. Let $X,Y \subseteq V(G)$ be disjoint such that $Y$ is stable and $|Y| \geq 2$ (possibly, $X = \emptyset$). Let $x,y \in V(G) \sm (X \cup Y)$. Suppose the following terms are satisfied:\\
\emph{\bfm{(\sref{rooted-nonseparating-path}.a)}} there is an induced $xy$-path in $G-(X \cup Y)$ with a bridge containing $Y$, and\\
\emph{\bfm{(\sref{rooted-nonseparating-path}.b)}} there is no induced $xy$-path in $G-(X \cup Y)$ with a bridge containing 
$Y \cup \{x\}$ for any $x \in X$.

Then, there is an induced $xy$-path $P$ in $G-(X \cup Y)$ such that:\\
\emph{\bfm{(\sref{rooted-nonseparating-path}.c)}} each $P$-bridge either meets $X$ or contains $Y$, and\\
\emph{\bfm{(\sref{rooted-nonseparating-path}.d)}} the interior of each leaf block of the $P$-bridge containing $Y$ meets $Y$, and\\
\emph{\bfm{(\sref{rooted-nonseparating-path}.e)}} if $|X|< \kappa(G) -2$, then each $P$-bridge meeting $X$ is a singleton.  
\end{statement}

\noindent
{\sl Proof.}
By \bfm{(\sref{rooted-nonseparating-path}.a)} and \bfm{(\sref{rooted-nonseparating-path}.b)}, there exist $4$-tuples $(P,B,H,\B)$ where $P$ is an induced $xy$-path with a $P$-bridge $B$ containing $Y$, $H \subseteq B$ contains $Y$ with each of its leaf blocks meeting $Y$ internally, and $\B$ is the union of $P$-bridges meeting $X$. Existence of $H$ follows from the fact that given a subgraph of $B$ containing $Y$ one may trim the block tree of such a subgraph by successively removing the interior of leaf blocks not meeting $Y$ internally. 

Choose a $4$-tuple $\X=(P,B,H,\B)$ such that\\
(I) $H$ is maximal,\\  
(II) subject to (I) $\B$ is maximal, and\\ 
(III) subject to (I) and (II) the number of $P$-bridges not meeting $X \cup Y$ is minimum.

Claims (\sref{rooted-nonseparating-path}.A,B,C)(see below) imply \sref{rooted-nonseparating-path}. \\

\noindent
(\sref{rooted-nonseparating-path}.A) \emph{In $\X$: there are no $P$-bridges not meeting $X \cup Y$.}\\
{\sl Proof.} Assume, towards contradiction, that there is a $P$-bridge $D$ not meeting $X \cup Y$ and let $v_{_1}$, $v_{_2}$ be its extremal attachments on $P$.
Let $C \not= D$ be a $P$-bridge adjacent to $(v_{_1}Pv_{_2})$; existence of which follows from $5$-connectivity. 

For an induced $v_{_1} v_{_2}$-path $P'$ in $G[V(D) \cup \{v_{_1}, v_{_2}\}]$, let $\D$ denote the union of the $P'$-bridges of $D$ not adjacent to $(v_{_1}Pv_{_2})$. Let $\D'$ denote the remaining $P'$-bridges of $D$. Choose $P'$ so that $\D$ is minimal and let $P''$ be the induced $xy$-path obtained from $P$ by replacing $(v_{_1}Pv_{_2})$ with $P'$.   

$\D$ is nonempty; for otherwise $(D-P'')\cup C \cup (v_{_1}Pv_{_2})$ is contained is some $P''$-bridge contradicting (III). Let then $D_{_1},\ldots,D_{_k}$ be the components of $\D$, and let $a_{_i}$, $b_{_i}$ be the extremal attachments of $D_{_i}$ on $P''$, for $1 \leq i \leq k$. Let $c_{_i} \in (a_{_i}P''b_{_i})$ adjacent to
$D - (P'' \cup D_{_i})$ or $(v_{_1}Pv_{_2})$; existence of $c_{_i}$ follows from $5$-connectivity for otherwise $\{a_{_i},b_{_i},v_{_1},v_{_2}\}$ is a $4$-disconnector of $G$. 

The vertex $c_{_i}$ is not adjacent to a component of $\D'$ or to $(v_{_1}Pv_{_2})$ for then the induced $v_{_1} v_{_2}$-path in $G[V(D) \cup \{v_{_1}, v_{_2}\}]$ obtained from $P'$ by replacing $a_{_i}P'b_{_i}$ with an induced $a_{_i}b_{_i}$-path in $G[V(D_{_i}) \cup \{a_{_i},b_{_i}\}]$ contradicts the choice of $P'$. Consequently, since $P''$ is induced, $N((a_{_i}P''b_{_i})) \subseteq V(\D)$, for $1 \leq i \leq k$. Therefore, $\bigcup_{i=1}^k (a_{_i}P''b_{_i})$ is a subpath of $P''$; let $a$ and $b$ be its ends. Since $\{a,b,v_{_1},v_{_2}\}$ is not a $4$-disconnector of $G$, there exists a $c \in (aP''b)$ adjacent to a component of $\D'$ or to $(v_{_1}Pv_{_2})$. Consequently, there is an $1 \leq i \leq k $ such that $c \in (a_{_i}P''b_{_i})$ which is a contradiction. \inQED \\

\noindent
(\sref{rooted-nonseparating-path}.B) \emph{In $\X$: $H$ coincides with $B$.}\\
{\sl Proof.} The $H$-bridges of $G-P$ are not adjacent to and do not meet members of $\B$, by \bfm{(\sref{rooted-nonseparating-path}.b)}. Consequently, the $H$-bridges of $G-P$ are subgraphs of $B$, by (\sref{rooted-nonseparating-path}.A). There are no $H$-bridges of $B$ with $\geq 2$ attachments in $H$ for the union of $H$ and an $H$-ear in $B$ linking two such attachments contradicts (I) as the property of each leaf block of such a union meeting $Y$ is sustained. 

It remains to prove that there are no $H$-bridges of $B$ with a single attachment in $H$. For suppose $D$ is such an $H$-bridge, let $v_{_1}$, $v_{_2}$ be its extremal attachments on $P$ (these exist due to $5$-connectivity), and let $v$ be its sole attachment in $H$. Suffices now that we prove the following claims: \\
(\sref{rooted-nonseparating-path}.B.1) \emph{No member of $\B$ is adjacent to $(v_{_1}Pv_{_2})$.}\\
(\sref{rooted-nonseparating-path}.B.2) \emph{Let $\P$ be the set of $((v_{_1}P v_{_2}),H)$-paths internally-disjoint of $D \cup P$. Then,}

(\sref{rooted-nonseparating-path}.B.2.1) \emph{$\P$ is nonempty; }

(\sref{rooted-nonseparating-path}.B.2.2) \emph{no member of $\P$ is adjacent to or meets a member of $\B$;} 

(\sref{rooted-nonseparating-path}.B.2.3) \emph{the members of $\P$ have a common end in $H$; and,}

(\sref{rooted-nonseparating-path}.B.2.4) \emph{there exists a member of $\P$ of order $\geq 3$.}\\

Indeed, assuming (\sref{rooted-nonseparating-path}.B.1-2) (see proofs below), one may argue as follows. 
Let $u \in V(H)$ be the common end of the members of $\P$, existence of which follows from (\sref{rooted-nonseparating-path}.B.2.3). By (\sref{rooted-nonseparating-path}.B.2.2) and (\sref{rooted-nonseparating-path}.B.2.4), there is an $H$-bridge $C$ of $B$ whose sole attachment in $H$ is $u$; and $C$ meets vertices of $\P$. Let $u_{_1}$ and $u_{_2}$ be its extremal attachments on $P$. Clearly, (\sref{rooted-nonseparating-path}.B.1-2) apply to $C$. 

Since $C$ meets $\P$, $u_{_1} P u_{_2} \cap v_{_1} P v_{_2} \not= \emptyset$. 
If $v_{_1} P v_{_2} \subseteq u_{_1} P u_{_2}$, then there exists a $((u_{_1} P u_{_2}),H)$-path internally-disjoint of $D \cup C \cup P \cup \{u,v\}$ with an end disjoint of $\{u,v\}$. For otherwise since $|V(H)| \geq 3$ (as $Y \subseteq V(H)$ and stable) and by (\sref{rooted-nonseparating-path}.B.1), $\{u,v,u_{_1},u_{_2}\}$ is a $4$-disconnector of $G$. Existence of such a path contradicts (\sref{rooted-nonseparating-path}.B.2.3). 

Consequently, $v_{_1} P v_{_2} \not\subseteq u_{_1} P u_{_2}$ and (by symmetry) $u_{_1} P u_{_2} \not\subseteq v_{_1} P v_{_2}$. Assume then, without loss of generality, that $x, u_{_1}, v_{_1}, u_{_2}, v_{_2}, y$ appear in this order along $P$. In an analogous manner to the previous case, by (\sref{rooted-nonseparating-path}.B.1), since $|V(H)| \geq 3$, and since $\{u,v,u_{_1},v_{_2}\}$ is not a $4$-disconnector of $G$, there is a $((u_{_1} P u_{_2}),H)$-path internally-disjoint of $D \cup C \cup P \cup \{u,v\}$ with an end disjoint of $\{u,v\}$. Since $v_{_1} P v_{_2}$ and $u_{_1} P u_{_2}$ overlap on $P$; existence of such a path contradicts (\sref{rooted-nonseparating-path}.B.2.3).

It remains to prove (\sref{rooted-nonseparating-path}.B.1-2). To see (\sref{rooted-nonseparating-path}.B.1), let $P'$ be an induced $xy$-path obtained from $P$ by replacing $v_{_1} P v_{_2}$ with an induced $v_{_1}v_{_2}$-path in $G[V(D) \cup \{v_{_1},v_{_2}\}]$. If a member of $\B$ is adjacent to $v_{_1} P v_{_2}$, then $P'$ contradicts (II). 

By (\sref{rooted-nonseparating-path}.B.1), since $|V(H)| \geq 3$, and since $\{v,v_{_1},v_{_2}\}$ is not a $3$-disconnector of $G$, $\P$ is nonempty. A member of $\P$ adjacent to or meeting a member of $\B$ implies that $P'$ (as defined above) contradicts \bfm{(\sref{rooted-nonseparating-path}.b)}. Next, if $Q,Q' \in \P$ do not have a common end in $H$, then the union of $Q$, $Q'$, and $v_{_1} P v_{_2}$ contains an $H$-ear. The path $P'$ and the union of $H$ and such an ear contradict (I).
Finally, if all members of $\P$ are edges, then by (\sref{rooted-nonseparating-path}.B.1), (\sref{rooted-nonseparating-path}.B.2.3), and since $|V(H)| \geq 3$, the set consisting of $v,v_{_1},v_{_2}$, and the common end of members of $\P$ is a $4$-disconnector of $G$. \inQED \\

\noindent
(\sref{rooted-nonseparating-path}.C) {\em $\X$ satisfies \emph{\bfm{(\sref{rooted-nonseparating-path}.e)}}.}\\
{\sl Proof.} Let $B' \in \B$ and let $v \in V(B) \cap X$. If $V(B') \setminus \{v\} \not= \emptyset$, then the union of components of $B'-v$, namely $\C(B',v)$ is nonempty. Let $C \in \C(B',v)$ and let $v_{_1},v_{_2}$ be the extremal attachments of $C$ on $P$; such vertices exist as $v \in N_G(C) \subseteq V(P) \cup \{v\}$ and $|N_G(C)| \geq 5$. 
Let $P'$ be an induced $xy$-path obtained from $P$ by replacing $v_{_1} P v_{_2}$ with an induced $v_{_1}v_{_2}$-path in $G[V(C) \cup \{v_{_1},v_{_2}\}]$.
\begin{equation}\label{H}
\mbox{$H$ (which coincides with $B$, by (\sref{rooted-nonseparating-path}.B)), is not adjacent to $(v_{_1} P v_{_2})$.}
\end{equation}
{\sl Proof.} Indeed, suppose to the contrary that $H$ is incident to $(v_{_1} P v_{_2})$.  Consider the following argument.

\begin{itemize}
\item[\bfm{(a)}] No member of $\C(B',v) \cup \B$ is adjacent to $(v_{_1} P v_{_2})$; for then $P'$ is a path contradicting \bfm{(\sref{rooted-nonseparating-path}.b)}. 

\item[\bfm{(b)}] Since $\{v_{_1},v_{_2},v\}$ is not a $3$-disconnector of $G$, the set $\P$ of $((v_{_1} P v_{_2}),H)$-paths internally-disjoint of $C \cup P$ is nonempty, by \bfm{(a)}. 

\item[\bfm{(c)}] Clearly, no member of $\P$ is adjacent to or meets a member of $\C(B',v) \cup \B$, for then $P'$ is a path contradicting \bfm{(\sref{rooted-nonseparating-path}.b)}. 

\item[\bfm{(d)}] Moreover, the members of $\P$ have a common end, say $u$, in $H$. Indeed, if not so and $Q,Q' \in \P$ have distinct ends in $H$, then the union of $Q$, $Q'$, and $(v_{_1} P v_{_2})$ contain an $H$-ear. The path $P'$ and the union of $H$ and such an ear contradict (I).

\item[\bfm{(e)}] Since $H$ coincides with $B$, the members of $\P$ are edges. It follows then that $\{v,u,v_{_1},v_{_2}\}$ is a $4$-disconnector of $G$. \inQED
\end{itemize}

To conclude our proof of (\sref{rooted-nonseparating-path}.C), let $X' \subseteq X$ satisfying (i) $|X' \cap V(B')|=1$ for every $B' \in \B$, and (ii) subject to (i)   
there is a $v \in X'$ such that $C(B,v) \not= \emptyset$ for some $B \in \B$. 

Let $\C = \bigcup_{B' \in \B, v \in X' \cap V(B')} \C(B',v)$. Let $\A$ denote the subpaths of $P$ each a union of extremal segments of overlapping and nested members of $\C$. By (\ref{H}), $H$ is not adjacent to $int A$ for any $A \in \A$. Let $A \in \A$ be maximal. The set comprised of $X$ and the ends of $A$ separates $H$ from the members of $\C$ defining $A$. As, by assumption, $|X| < \kappa(G)-2$, such a set has size $< \kappa(G)$. \inQED \\
\QED

As indicated above, in this paper, the counterpart of \sref{rooted-nonseparating-path} is the following.

\begin{statement}\label{always}
Suppose $\kappa(G)\geq 5$ and  $H\iso K_{_{2,3}}\subseteq G$ such that the $3$-valent vertices of $H$ are not adjacent in $G$. Then $G$ contains an induced path $P$ meeting no $2$-valent vertex of $H$ whose ends are the $3$-valent vertices of $H$ and with a $P$-bridge meeting at least two $2$-valent vertices of $H$. 
\end{statement}

\proof
Let $V(H) = \{x_{_i}\}_{i\in [5]}$ such that $d_H(x_{_1}) = d_H(x_{_2}) = 3$. 
Existence of disjoint $x_{_1} x_{_2}$-path and an $x_{_4} x_{_5}$-path in $G-x_{_3}$, clearly implies the assertion. By the main result of \cite{seymour,thomassen} (the characterization of $2$-linked graphs), nonexistence of such paths in $G-x_{_3}$ occurs if and only if $G-x_{_3}$ has an embedding in the plane such that $x_{_1}, x_{_4}, x_{_2}, x_{_5}$ appear on the outerface, say $X$ (which is a circuit), in this clockwise order. As $\kappa(G-x_{_3}) \geq 4$ and $x_{_1}, x_{_4}, x_{_2}, x_{_5}$ induce a $C_{_4}$, $X \iso C_{_4}$. 

Let $C$ be an $\{x_{_1},x_{_2}\}$-circuit in $G-\{x_{_3},x_{_4},x_{_5}\}$, and let $P_{_1},P_{_2}$ be the two $x_{_1}x_{_2}$-paths comprising $C$.
Since $x_{_1} x_{_2} \notin E(G)$, $int P_{_i} \not= \emptyset$, $i =1,2$.  
The circuits $C$ and $X$ divide the plane into $4$ regions:\\
$R_{_1} \colon =$ $int C_{_1}$ with $C_{_1}$ comprised of, say, $P_{_1}$ and the edges $\{x_{_5}x_{_1},x_{_5}x_{_2}\}$.\\
$R_{_2}\colon =$ $int C$ (not including $C$).\\
$R_{_3}\colon=$ $int C_{_2}$ with $C_{_2}$ comprised of, say, $P_{_2}$ and the edges $\{x_{_4}x_{_1},x_{_4}x_{_2}\}$.\\
$R_{_4} \colon = $ $ext X$.

Let $y \in N_G(x_{_3})\sm \{x_{_1},x_{_2}\}$ embedded in $R_{_i}$, $i \in =1,2,3$. Since $\kappa(G-x_{_3}) \geq 4$, there is a $(y,int P_{_i})$-path, $Q$, not meeting $P_{_{3-i}}$ for some $i =1,2$. Let $x \in V(C_{_i}) \cap \{x_{_4},x_{_5}\}$. There is an $(x,int P_{_i})$-path, $Q'$, not meeting $P_{_{3-i}}$. Since $Q\cup Q' \cup P_{_i} \cup \{x_{_3}\}$ contains an $x x_{_3}$-path disjoint of $P_{_{3-i}}$ the claim follows. \QED

\subsection{Salami and Pie}\label{salami-and-pie}
The purpose of this section is to state the main result of Yu \cite{Yu} that provides a characterization of the $4$-connected plane graphs containing no $TK_{_4}$ rooted at a prescribed set of $4$ vertices. 

Unless otherwise stated, throughout this section $G$ is a $4$-connected plane graph, and $W =\{w_{_1},w_{_2},w_{_3},w_{_4}\} \subseteq V(G)$.\\

\noindent
\scaps{Definition B.}\\ A connected plane graph $J$ is called a \emph{configuration}
if it has $6$ distinct vertices $\{a_{_1},a_{_2},a_{_3},a_{_4},x,y\}$ such that:\\
(B.1) $A=\{a_{_1},a_{_2},a_{_3},a_{_4}\}$ is stable and appears in this clockwise order on the outerwalk of $J$,\\
(B.2) no disconnector $T$ of $J$, $|T| \leq 3$, separates $x$ or $y$ from $A \sm T$,\\
(B.3) no disconnector $T$ of $J$, $ |T| \leq 4$, separates $\{x,y\}$ from $A \sm T$, and\\
(B.4) every disconnector $T$ of $J$, $|T| \leq 3$, separates two vertices of $A \cup \{x,y\}$.\\

We refer to $x$ and $y$ as the \emph{roots} of $J$ and to the members of $A$ as its \emph{exits}. 
  
\begin{figure}[htbp]
	\centering
		\includegraphics{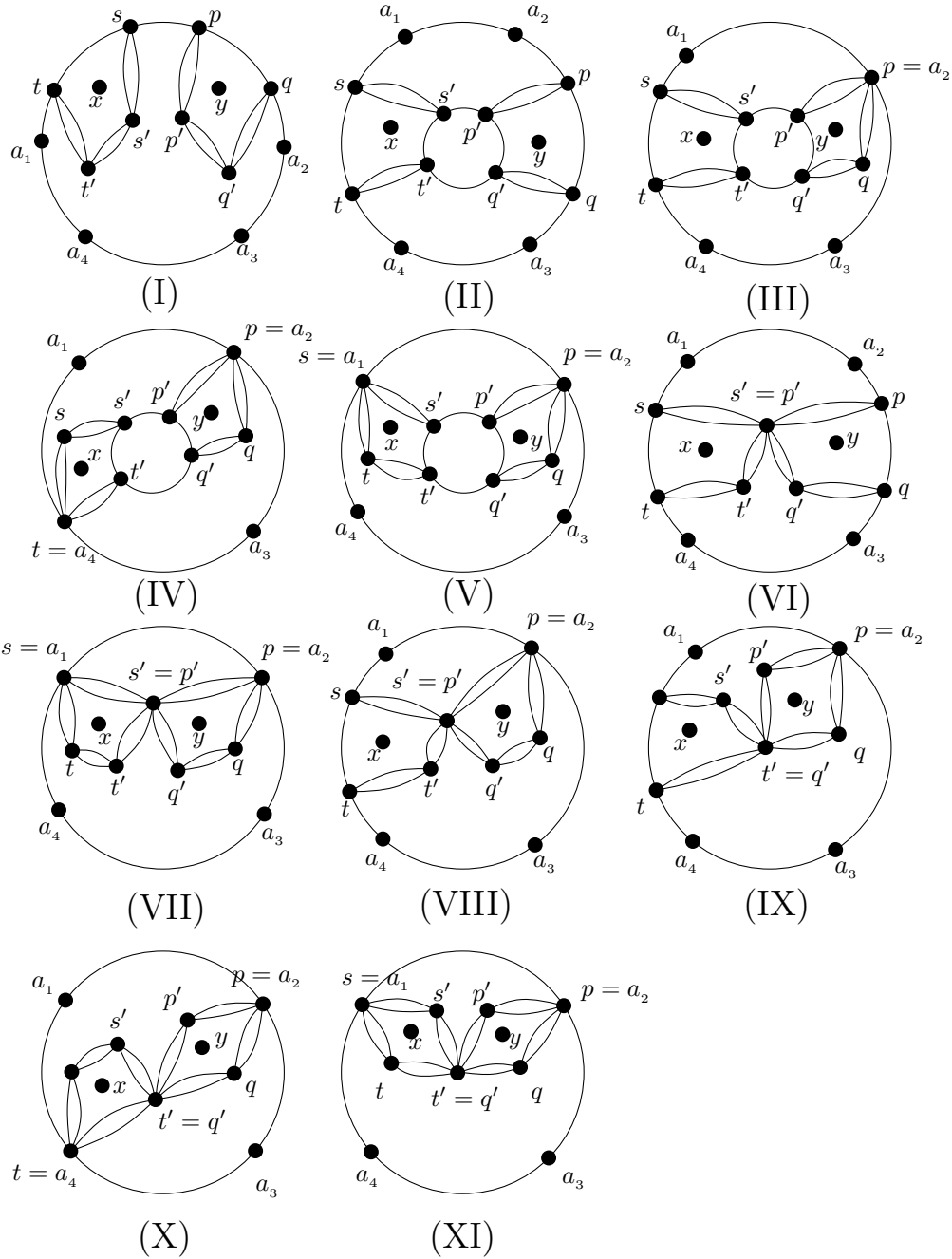}
	\caption{Yu's configurations (copy of \cite[Fig.~2]{Yu}).}
	\label{confs}
\end{figure}

Each graph in Figure~\ref{confs} is a configuration with outerwalk, say $C$. 
The configurations (II)-(XI) each contains a $4$-disconnector $S_{_x}=\{s,s',t,t'\}$ separating $(A \sm S_{_x})\cup \{y\}$ from $x$, and a $4$-disconnector $S_{_y}=\{p,p',q,q'\}$ separating $(A \sm S_{_y})\cup \{x\}$ from $y$. The component of $J-S_{_x}$ conataining $x$ is disjoint of the component of $J-S_{_y}$ containing $y$. 
In configuration (I), if $x \notin V(C)$, then $S_{_x}$ as above exists. Otherwise, if $x \in V(C)$, we put $S_{_x} =\{x\}$. A similar convention holds for $y$. We assume that $\{t,s,s',t'\}$ and $\{p,q,q',p'\}$ appear around $x$ and $y$ in these clockwise orders, respectively, as depicted in Figure~\ref{confs}. 

In configuration (I), $a_{_1},t,s,p,q,a_{_2},a_{_3},a_{_4}$ appear in this clockwise order on $C$; and possibly $a_{_1}=t$ or $s=p$ or $q= a_{_2}$. In configuration (II), $a_{_1},a_{_2},p,q,a_{_3},a_{_4},t,s$ appear in this clockwise order on $C$; and possibly $s= a_{_1}$ or $t=a_{_4}$ or $p=a_{_2}$ or $q = a_{_3}$. In addition, $s',p',q',t'$ lie on another face (other than $C$) and appear in this clockwise order on this face; possibly $s'=p'$ or $t'=q'$. 

Configurations (III)-(XI) are variations of (II). 
In these configurations some constraints of (II) are relaxed. 
For instance, the presence of the face (other than $C$) containing $s',p',q',t'$
is substituted by the constraint that $s'=p'$ or $t'=q'$. Also,
it is no longer required that $q$ and $t$ lie on $C$ and so on. 
We refer and advise the reader to consult \cite{Yu} for complete details.  

It will be useful for us to observe the following.

\begin{statement}\label{common1}
If $J$ is a configuration as in Figure~\ref{confs}, then no exit of $J$ is a common neighbor of its roots. 
\end{statement}

\proof
An exit of $J$ serving as a common neighbor of the roots of $J$ contradicts (B.2).  
Let $x,y, A= \{a_{_i}\}_{i=1,2,3,4}, S_{_x},S_{_y}$, and $C$ be as above. 
Let $a \in N(x) \cap N(y) \cap A$.

If $S_{_x} \not=\{x\}$ and $S_{_y}\not=\{y\}$, then $a \in A \cap S_{_x} \cap S_{_y} \not= \emptyset$; implying that $s'=p' \in A$ or $t'=q' \in A$. In each configuration, the latter translates into a disconnector $T$, $|T| \leq 3$, of $J$ separating at least one of $x,y$ from $A\sm T \not= \emptyset$; contradicting (B.2). 

Suppose then that $S_{_x} = \{x\}$; implying that $J$ is configuration $(I)$ and $x \in V(C)$.
If $y \in V(C)$ as well, then at least one of the edges $xa,ya$ forms a $2$-disconnector $T$ of $J$ separating $x$ or $y$ from $A\sm T \not= \emptyset$. Assume then that $S_{_y} \not= \emptyset$. Consequently, $a_{_1} = a$, for otherwise $\{x,a\}$ is a $2$-disconnector separating $x_{_1}$ from $y$. Moreover, $a \in \{p',q'\}$. The latter implies that there is a disconnector $T \subseteq \{p',q',q\}$, separating  $y$  (and $x$) from $A \sm T$; contradicting (B.2).\QED \\

\noindent
\scaps{Definition C.}\\ A pair $(G,W)$ is called a \emph{pie} if there is a facial circuit $C$ of $G$, called the \emph{mold}, and sets $S_{w_{_i}}$, $i =1,2,3,4$, such that either:\\
(C.1) $w_{_i} \in V(C)$ and then $S_{w_{_i}} = \{w_{_i}\}$, or\\
(C.2) $G$ has a $4$-disconnector $S_{w_{_i}}$ satisfying $|S_{w_{_i}} \cap V(C)| =2$, and separating $w_{_i}$ from $W \sm \{w_{_i}\}$. The components of $G-S_{w_{_i}}$ containing $w_{_i}$ are disjoint and $w_{_j} \notin S_{w_{_i}}$ whenever $i \not= j$.\\

\begin{figure}[htbp]
	\centering
		\scalebox{0.7}{\includegraphics{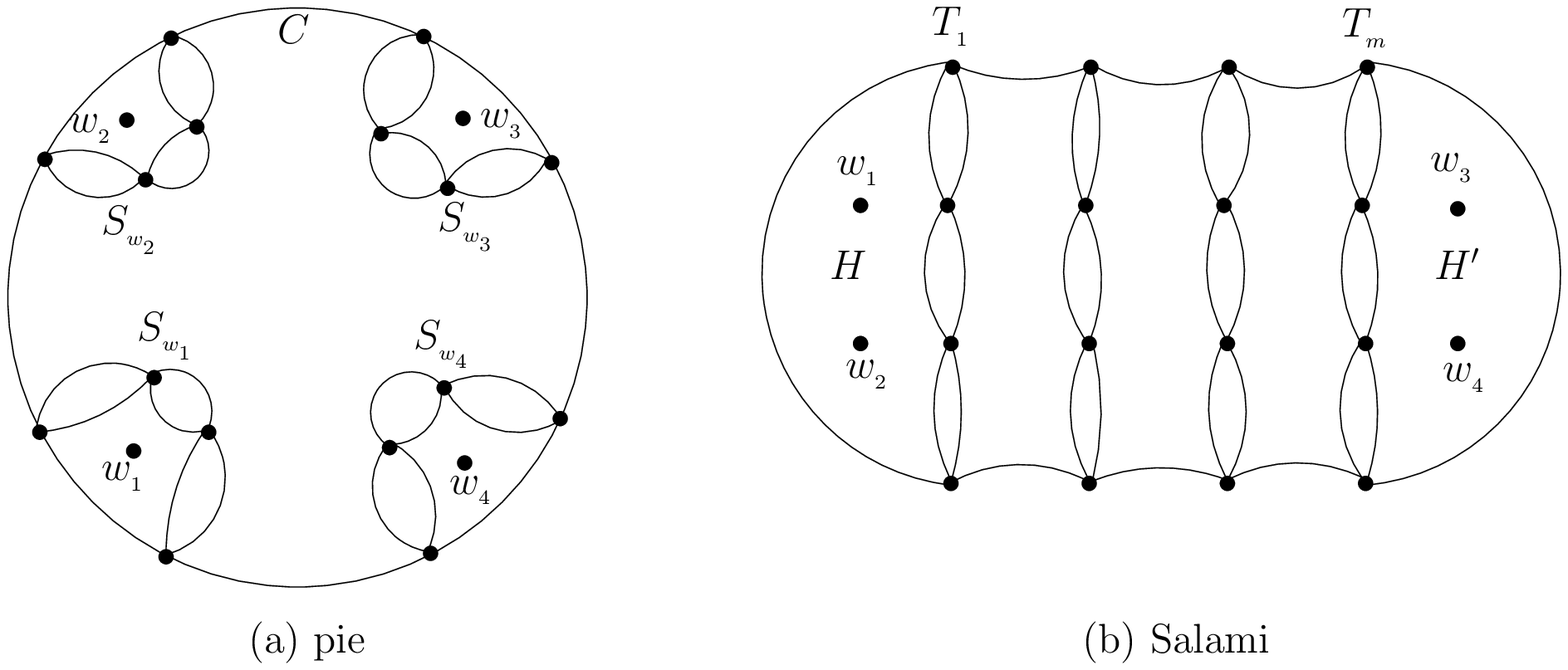}}
	\caption{Yu's Salami and pie (copy of \cite[Fig.~7]{Yu}).}
	\label{fig:salami-pie}
\end{figure}

An illustration of a pie can be found in Figure~\ref{fig:salami-pie}(a). 
In the figure, the mold is drawn as the outerface. \\

\noindent
\scaps{Definition D.}\\ A pair $(G,W)$ is called a \emph{Salami separating $\{w_{_1},w_{_2}\}$ from $\{w_{_3},w_{_4}\}$} if $G$ contains 
distinct (not necessarily disjoint) $4$-disconnectors $T_{_i}$, $i \in [m]$
such that:\\
(D.1) each $T_{_i}$ separates $\{w_{_1},w_{_2}\}$ from $\{w_{_3},w_{_4}\}$,\\
(D.2) each $T_{_i}$ separates $\{w_{_1},w_{_2}\}$ from $T_{_{i+1}}\sm T_{_i}\not= \emptyset$,\\
(D.3) $G$ has no $4$-disconnector $T$ separating $T_{_i} \sm T \not= \emptyset$ from $T_{_{i+1}} \sm T \not= \emptyset$,\\
(D.4) the $4$-hammocks $H$ and $H'$ satisfying $\bnd{H}{G}=T_{_1}$ and $\bnd{H'}{G} = T_{_m}$ are minimal subject to containing $\{w_{_1},w_{_2}\}$ and $\{w_{_3},w_{_4}\}$, respectively. \\
(D.5) $J \in \{H,H'\}$ is a configuration as listed in Figure~\ref{confs} such that there are bijections $\bnd{J}{G} \rightarrow A$ and 
$W \cap V(J) \rightarrow \{x,y\}$, where $A$ and $x,y$ are as in \scaps{Definition B}.\\

An illustration of a Salami can be found in Figure~\ref{fig:salami-pie}(b). We refer to the hammocks $H,H'$ of (D.4) as the \emph{ends} of the Salami $(G,W)$.

A pair $(G,W)$, whether a pie or a Salami, has a set $S_{w}$ associated with each $w \in W$. 
Such a set, if satisfies $S_{w} = \{w\}$, is called \emph{elementary}; otherwise, $S_{w}$ is the boundary of a $4$-hammock containing $w$ in its interior, and is then called \emph{nonelementary}. A  $4$-hammock containing $w$ and bounded by a nonelementary $S_{_w}$ has the members of $S_{_w}$ appear in a clockwise order along its outerwalk. 
We refer to this order as the \emph{order} of $S_{_w}$. 

\begin{statement}\label{rooted-tk4} 
\bfm{(Yu \cite[Theorem~$\mbs{4.2}$]{Yu})}\\
$G$ contains no $TK_{_4}$ rooted at $W$ if and only if $(G,W)$ is a Salami or a pie. 
\end{statement}    

In fact, Yu's argument for \sref{rooted-tk4}~\cite[Theorem 4.2]{Yu}, which relies on 
~\cite[Theorems 2.1, 3.1]{Yu}, establishes the following more precise and useful fact.

\begin{statement}\label{actually}
If $G$ contains no $TK_{_4}$ rooted at $W$ and there is a $4$-disconnector $T$ separating $\{w_{_1},w_{_2}\}$ from $\{w_{_3},w_{_4}\}$, then\\
\indent \emph{(\sref{actually}.A)} $(G,W)$ is a Salami separating
$\{w_{_1},w_{_2}\}$ from $\{w_{_3},w_{_4}\}$, and\\
\indent \emph{(\sref{actually}.B)} the ends of $(G,W)$ are contained in distinct $4$-hammocks whose boundary is $T$.
\end{statement}

A consequence of \sref{common1} that will be useful for us is the following. 

\begin{statement}\label{common2}
If three members of $W$ have a common neighbor not in $W$, then $(G,W)$ has a $TK_{_4}$ rooted at $W$. 
\end{statement}

\proof
Suppose $c \in (\bigcap_{i=1}^3 N(w_{_i})) \sm W$. 
By \sref{common1}, $(G,W)$ is not a Salami. Indeed, since any $4$-disconnector separating two pairs of $W$ contains $c$, both ends of $(G,W)$ contain $c$ as an exit; implying that at least one end has an exit that is a common neighbor of its roots, contradicting \sref{common1}. 

Suppose then that $(G,W)$ is a pie and let $C$ be its mold. Suppose $\{S_{w_{_1}} \cap V(C), S_{w_{_2}} \cap V(C),S_{w_{_3}} \cap V(C),S_{w_{_4}} \cap V(C)\}$ appear in this clockwise order on $C$. Since $\{w_{_1},c,w_{_3}\}$ is not a $3$-disconnector of $G$, at least one of $w_{_1}, w_{_3}$ is not on $C$; so let $S_{w_{_1}}$ be nonelementary. 

If $w_{_3} \notin V(C)$, then $c \in S_{w_{_1}} \cap S_{w_{_3}}$ and then $G$ has a disconnector $T \subset (S_{w_{_1}} \cap V(C)) \cup (S_{w_{_3}} \cap V(C)) \cup \{c\}$, $|T| \leq 3$, separating $w_{_2}$ from $w_{_4}$. On the other hand, if $w_{_3} \in V(C)$ then $G$ has a disconnector $T \subset (S_{w_{_1}} \cap V(C)) \cup \{c,w_{_3}\}$, $|T| \leq 3$, separating $w_{_2}$ from $w_{_4}$.\QED \\

We shall use the following observation that is a consequence of \sref{common2}.

\begin{statement}\label{notapex}
If $H \iso K_{_{2,3}} \subseteq G$ and $G$ is a nonplanar $5$-connected apex graph, then either $TK_{_5} \subseteq G$, or no $3$-valent vertex of $H$ is an apex vertex of $G$.  
\end{statement}

\proof
Fix an embedding of $G$ in the plane and treat $G$ as a plane graph. Let $u \in V(H)$ be $3$-valent and let $W \subset N_G(u)$ consist of the three $2$-valent vertices of $H$ and a fourth vertex not in $V(H)$. Since $3$ members of $W$ have a common neighbor not in $W$, it follows, by \sref{common2}, that if $u$ is an apex vertex of $G$, then $G-u$ contains a $TK_{_4}$ rooted at $W$; implying that $TK_{_5} \subseteq G$.
\QED \\

The facial circuits of a $4$-connected planar graph are its induced nonseparating circuits and are unique \cite{mt}. Let $C \subset G$ be an induced $4$-circuit with $V(C) =\{x_{_1},x_{_2},x_{_3},x_{_4}\}$ such that $x_{_1} x_{_3} \notin E(G)$. Let $H$ be a $4$-hammock of $G$ such that $x_{_1} \in int_G H$ and $x_{_3} \notin V(H)$. Clearly, $x_{_2},x_{_4} \in \bnd{H}{G} =\{u_{_0},u_{_1},u_{_2},u_{_3}\}$, where the latter appear in this clockwise order along the outerwalk of $H$ (here we refer to the embedding of $H$ induced by the embedding of $G$). It is not necessarily true that $C$ is a facial circuit of $G$. However, 
\begin{equation}\label{face}
\mbox{ $\{x_{_2},x_{_4}\} = \{u_{_{i(mod\; 4)}},u_{_{i+1(mod\; 4)}}\}$, $1 \leq i \leq 4$, if and only if $C$ is a facial circuit of $G$.} 
\end{equation}
To see (\ref{face}), suppose that $\{x_{_2},x_{_4}\} = \{u_{_{i(mod\; 4)}},u_{_{i+1(mod\; 4)}}\}$, $i =\{1,2,3,4\}$ and $C$ is a separating
circuit of $G$. That is, at least one  of $H \cap int C$, $(G-H) \cap int C$ is a nonempty graph. By the assumption on $\{x_{_2},x_{_4}\}$ and since $H$ is a hammock, at least one of $\{x_{_1},x_{_2},x_{_4}\}$, $\{x_{_2},x_{_3},x_{_4}\}$ is a $3$-disconnector of $G$; a contradiction. The converse is clearly true. 

Suppose, next, that $J$ is a configuration as in Figure~\ref{confs} with $x_{_1}$ and $x_{_3}$ as its roots such that $x = x_{_1}$ and $y=x_{_3}$, say, and $S_{_x},S_{_y}$ are both nonelementary. Clearly, $x_{_2},x_{_4} \in S_{_x} \cap S_{_y}$. If (B.2-4) are to be satisfied subject to the allowed vertex identifications specified for configurations in Figure~\ref{confs}, then\\
(M.1) $\{x_{_2},x_{_4}\} = \{s,s'\} = \{p,p'\}$ in configuration (I), or\\
(M.2) $\{x_{_2},x_{_4}\} = \{s',t'\} = \{p',q'\}$ in configurations (II)-(XI).

A similar situation arises when only one of $S_{_x}$ and $S_{_y}$ is nonelementary and also when $S_{_x}$ and $S_{_y}$ are part of a pie and not a configuration. 
Consequently, we observe the following.

\begin{statement}\label{conf-face}
Let $x,y$ be part of an induced $4$-circuit such that $x$ and $y$ are nonadjacent. 
Let $z \in \{x,y\}$. 

\emph{(\sref{conf-face}.A)} If $J$ is a configuration as in Figure~\ref{confs} with roots $x,y$, then an  
$S_{_z}$ that is\\ \indent \indent \indent \ \ \ nonelementary satisfies $V(C) \sm \{x,y\} \subset S_{_z}$ and the members of $V(C) \cap S_{_z}$ appear\\ \indent \indent \indent \ \ \ as in (M.1) and (M.2).

\emph{(\sref{conf-face}.B)} If $(G,W)$ is a pie with mold $F$, where $W = \{x,y,\ell,\ell'\}$ and $|V(C) \cap W|=2$, \indent \indent \indent \ \ \ then $\{S_{_x} \cap F, S_{_y} \cap F, S_{_{\ell}} \cap F, S_{_{\ell'}} \cap F\}$ appear in this clockwise order along $F$; an \indent \indent \indent \ \ \ $S_{_z}$ that is nonelementary satisfies $V(C) \sm \{x,y\} \subset S_{_z}$ such that the members of \indent \indent \indent \ \ \ $V(C) \cap S_{_z}$ appear consecutively in the order of $S_{_z}$: one on $F$ and the other not \indent \indent \indent \ \ \ on $F$. 
\end{statement}

The following is now an exercise. 

\begin{statement}\label{must-be-a-face}
Suppose $w_{_1},w_{_2}$ are a part of an induced $4$-circuit $C$ of $G$ such that 
$w_{_1} w_{_2} \notin E(G)$ and $|V(C) \cap W|=2$. Then, $G$ has a $TK_{_4}$ rooted at $W$ or $C$ is a facial circuit of $G$. 
\end{statement}

\subsection{Proof of \sref{1con-1}}
Throughout this section, $G$ is a $5$-connected nonplanar apex graph satisfying $TK_{_5} \not\subseteq G$; $K$ is as in (\ref{K}); $a$ is an apex vertex of $G$ as in \sref{1con-1}; and  $P$ is an induced $x_{_1}x_{_2}$-path satisfying (\ref{P}).

By \sref{MaYuk4} and \sref{disc}, to prove \sref{1con-1}, suffices that we prove the following claims:\\
(\sref{1con-1}.A) \emph{If $a \in \{x_{_3},x_{_4},x_{_5}\}$, then $K^-_{_4} \subseteq G$ or $G$ has a planar hammock or $TK_{_5} \subseteq G$.}\\ 
(\sref{1con-1}.B) \emph{If $a \notin V(K)$, then $K^-_{_4} \subseteq G$ or $G-a$ is nonplanar.}\\

Claim (\sref{1con-1}.B) is an exercise; though not necessarily a short one. Nevertheless, it becomes a routine case analysis due to (\ref{P}) (which is a consequence of  \sref{rooted-nonseparating-path}). We concentrate on (\sref{1con-1}.A) only. 

\subsubsection{Proof of (\sref{1con-1}.A)}
Throughout this proof $a = x_{_3}$ and $C$ is the $C_{_4}$ induced by $V(K) \sm \{a\}$. The facial circuits of a $4$-connected planar graph are its induced nonseparating circuits and are unique \cite{mt}. Consequently, by (\ref{stable}) and \sref{must-be-a-face}, 
\begin{equation}\label{setting1}
\mbox{$G-P$ is a chain, $a$ is a cut vertex of $G-P$, and $C$ is a facial circuit of $G-a$.}
\end{equation}
To see (\ref{setting1}), note that if $G-P-a$ contains an $x_{_4} x_{_5}$-path $P'$, then, in $G-a$, the set $\{x_{_1},x_{2},x_{_4},x_{_5}\}$ forms an induced $4$-circuit separating  $int P$ from $int P'$ each is nonempty by (\ref{stable}). Therefore, by \sref{must-be-a-face} applied to $(G-a,\{x_{_1},x_{_2},y,y'\})$, $y,y' \in N_G(a) \sm V(K)$ (these exist by (\ref{stable})), it follows that $TK_{_5} \subseteq G$. This establishes that $G-P$ is a chain and that $a$ is its cut vertex. Consequently, $G-a$ has an embedding in which $P$ and $G-P-a$ are both embedded in $int C$ or $ext C$ and thus $C$ is a facial circuit of $G-a$.\\

\noindent
\scaps{Outline.} We now outline the reminder of our proof for (\sref{1con-1}.A); adjourning technical details until later parts of this section. Throughout, the remainder of this section $W$ is a set of the form $\{x_{_1},x_{_2},\ell,\ell'\}$, where $\ell,\ell' \in N_G(a) \sm V(K)$, $G-a$ is embedded in the plane such that $C$ is its outer face; such an embedding is called $\Pi$. We begin by proving (see proof below) that\\

\noindent
(\sref{1con-1}.A.1) \emph{there exist $\ell,\ell' \in N_G(a) \sm V(K)$ 
such that $(G-a,\{x_{_1},x_{_2},\ell,\ell'\})$ is a Salami.} \\

We distinguish between two cases:\\
\scaps{Case A.} Either $G-a$ has a $4$-disconnector $T$ separating $\{x_{_1},\ell\}$ from $\{x_{_2},\ell'\}$ for some $\ell,\ell' \in N_G(a) \sm V(K)$; or\\
\scaps{Case B.} $T$ as above does not exist in $G-a$ for any $\ell,\ell' \in N_G(a) \sm V(K)$.  

Our argument for \scaps{Case A} is as follows. Let $\L$ be the Salamies
$(G-a,W)$ separating $\{x_{_1},\ell\}$ from $\{x_{_2},\ell'\}$ for some $\ell,\ell' \in N_G(a) \sm V(K)$. The assumption of this case that $T$ exists and \sref{actually} imply that $\L \not= \emptyset$. For $(G-a,W) \in \L$, let $\alpha(G-a,W) = |V(J) \cup V(J')|$, where $J$ and $J'$ are the ends of $(G-a,W)$.  

\begin{equation}\label{choice}
\mbox{Choose $(G-a,W) \in \L$ minimizing $\alpha(G-a,W)$.}
\end{equation} 

\noindent
\scaps{Setting A.} For $(G-a,W)$ as in (\ref{choice}), let $J_{_i}$ denote the end of $(G-a,W)$ containing $x_{_i}$, $i=1,2$.
Put $x \in \{x_{_1},x_{_2}\}$ and $y \in \{\ell,\ell'\}$. Let $J \in \{J_{_1},J_{_2}\}$ such that $x,y \in int_{G-a} J$. Let $S_x$, $S_y$, $A = \{a_{_1},a_{_2},a_{_3},a_{_4}\} = \bnd{J}{G-a}$ be as these are defined in \S\ref{salami-and-pie} and Figure~\ref{confs}. 
Since any $4$-disconnector separating $x_{_1}$ and $x_{_2}$ contains $\{x_{_4},x_{_5}\}$,
we have that $\{x_{_4},x_{_5}\} \subseteq bnd_{G-a} J$.\\ 

A proof of (\sref{1con-1}.A) then follows from the following two claims (that are proven below): 

\noindent
(\sref{1con-1}.A.2) \emph{$J-\{x_{_4},x_{_5}\}$ contains two paths $Y$ and $Q$ such that:\\
\indent \emph{(\sref{1con-1}.A.2.1)} $Y$ is an $xa_{_2}$-path and $Q$ is an $xa_{_3}$-path; and,\\
\indent \emph{(\sref{1con-1}.A.2.2)} $V(Y) \cap V(Q) = \{x,y\}$.}

\noindent
(\sref{1con-1}.A.3) \emph{$J-\{x_{_4},x_{_5}\}$ contains an $(x,\{a_2,a_{_3}\})$-$2$-fan with one of its members containing $y$.}\\

\noindent
\bfm{Proof that} (\sref{1con-1}.A.2-3) \bfm{imply} (\sref{1con-1}.A)\bfm{.}   
Let $bnd_{G-a} J_{_1}  = U = \{u_{_1},u_{_2},u_{_3},u_{_4}\}$ and let $bnd_{G-a} J_{_2}=Z= \{z_{_1},z_{_2},z_{_3},z_{_4}\}$ such that $u_{_1} = x_{_4} = z_{_1}$, $u_{_4}=x_{_5}= z_{_4}$, $\{u_{_1},u_{_2},u_{_3},u_{_4}\}$ and $\{z_{_1},z_{_4},z_{_3},z_{_2}\}$ appear in these clockwise order on the outerwalks of $J_{_1}$ and $J_{_2}$, respectively, with respect to $\Pi$. 

Let $Y,Q \subseteq J_{_1}$ be as in (\sref{1con-1}.A.2), where $x =x_{_1}$,$y = \ell$, $a_{_2} = u_{_2}$, and $a_{_3} = u_{_3}$.
Let $F \subseteq J_{_2}$ be as in (\sref{1con-1}.A.3) where $x = x_{_2}$, $y= \ell'$, $a_{_2} =z_{_2}$, and $a_{_3} =z_{_3}$. Let $Y'$ be the member of $F$ containing $\ell$ such that, without loss of generality, $z_{_2} \in V(Y')$. Let $Q'$ be the other member of $F$. 
Let $F'=\{Y'',Q''\}$ be a $(\{u_{_2},u_{_3}\},\{z_{_2},z_{_3}\})$-$2$-linkage in $G-int_{G-a}J_{_1} -int_{G-a}J_{_2} - \{a,x_{_4},x_{_5}\}$. Such a linkage clearly exists as $\kappa(G-\{a,x_{_4},x_{_5}\}) \geq 2$. By planarity, we may assume that $Y''$ is a $u_{_2} z_{_2}$-path and that $Q''$ is a $u_{_3} z_{_3}$-path. 
Let $X_{_Y}$ and $X_{_Q}$ be the $x_{_1} x_{_2}$-paths in $G-a$ that are the union of $Y,Y',Y''$ and $Q,Q',Q''$, respectively.  
Let $F''$ be an $(x_{_4},X_{_Y})$-$4$-fan in $G-a$ with two of its members the edges $\{x_{_1}x_{_4},x_{_2}x_{_4}\}$ and minimizing the length of its members. Let $X,X'$ be the remaining members of $F''$. Such do not meet $X_{_Q}$, by planarity.

Now, the path $X_{_Q}$ satisfies (F.2); implying that $TK_{_5} \subseteq G$. 
To see this, suffices to show that $G-X_{_Q}-a+\ell$ contains an $\ell \ell'$-path containing $x_{_4}$ (recall that $x_{_3} \ell,x_{_3} \ell' \in E(G)$). Such a path clearly exists in $X \cup X' \cup (X_{_Y}-\{x_{_1},x_{_2}\})$.\QED \\ 

With \scaps{Case A} resolved, we may now assume, due to \sref{actually}, that in $G-a$ 
\begin{equation}\label{nocut}
\mbox{there are no $\ell,\ell' \in N_{G}(a)\sm V(K)$ and a $4$-disconnector separating $\{x_{_1} ,\ell\}$, $\{x_{_2},\ell'\}$.}
\end{equation} 
Consequently, by (\sref{1con-1}.A.1), there are $\ell,\ell' \in N_{G}(a)\sm V(K)$ such that $(G-a,W)$ is a Salami separating $\{x_{_1},x_{_2}\}$ from $\{\ell,\ell'\}$. This and (\ref{nocut}) clearly imply that\\

\noindent
(\sref{1con-1}.A.4) \emph{there exist $\ell,\ell' \in N_{G}(a) \sm V(K)$ such that $(G-a,W)$ is a Salami such that\\
\indent \emph{(\sref{1con-1}.A.4.1)} $x_{_1},x_{_2}$ are contained in a common end, and\\
\indent \emph{(\sref{1con-1}.A.4.2)} for $i=1,2$, $S_{x_{_i}}$, if nonelementary, defines a $4$-hammock $H_{_i}$ of $G-a$ satisfying
\begin{center} 
$x_{_i} \in int_{G-a} H_{_i}$, $\bnd{H_{_i}}{G-a} = S_{x_{_i}}$, and $N_G(a) \cap int_{G-a}H_{_i} = \{x_{_i}\}$.\end{center}}

\noindent
\bfm{Proof that} (\sref{1con-1}.A.4) \bfm{implies} (\sref{1con-1}.A)\bfm{.} 
Let $J$ be the end of $(G-a,W)$ containing $\{x_{_1},x_{_2}\}$. 
Clearly, a nonelementary $S_{x_{_i}}$ satisfies $\{x_{_4},x_{_5}\} \subset S_{x_{_i}}$. 
We prove the following claims \bfm{(a)} and \bfm{(b)}.\\

\noindent
\bfm{(a)} \emph{$S_{x_{_i}}$ is nonelementary for at least one $i=1,2$.}\\ 
{\sl Proof.} Otherwise,
the end $J$ of $(G-a,W)$ containing $x_{_1},x_{_2}$ is configuration (I) from Figure~\ref{confs} with both $x_{_1},x_{_2}$ appearing on its outerwalk. Since $x_{_1}x_{_2} \notin E(G)$, by (\ref{stable}), at least one of $x_{_4},x_{_5}$ also lies on the outerwalk 
of $J$ or $\{x_{_1},x_{_2},x_{_4},x_{_5}\}$ contains a $3$-disconnector of $G-a$. As
$\{x_{_1},x_{_2},x_{_4},x_{_5}\}$ define an induced $C_{_4}$, at least one of 
$\{x_{_4},x_{_5}\}$ has degree at most $2$ in $G-a$ which is a contradiction. \inQED\\

\noindent
\bfm{(b)} \emph{If $H_{_i}$ is as in \emph{(\sref{1con-1}.A.4.2)} for some $i=1,2$, then 
$H_{_i}$ is a trivial $5$-hammock of $G$ so that \indent $N_{G}(x_{_i}) = \bnd{H_{_i}}{G-a} \cup \{a\}$.}\\
{\sl Proof.} Since $N_G(a) \cap int_J H_{_i} = \{x_{_i}\}$, by (\sref{1con-1}.A.4.2), then $\bnd{H_{_i}}{G} = \bnd{H_{_i}}{G-a} \cup \{x_{_i}\} = S_{x_{_i}} \cup \{x_{_i}\}$. As $|V(G)| > |V(H_{_i})|$, we may assume that $H_{_i}$ is not a planar hammock or we are done; and thus, $|V(H_{_i})| \leq 6$. If $|V(H_{_i})|=6$, then the (sole) vertex $v$ in $int_G H_{_i}$ has $N_G(v) = \bnd{H_{_i}}{G}$; implying that $\{v,x_{_1},x_{_4},x_{_5}\}$ induce a $K^-_{_4}$.
As $H_{_i}$ is a $4$-hammock of $G-a$, it follows that $H_{_i}$ is a trivial $5$-hammock of $G$ and that $N_{G}(x_{_i}) = \bnd{H_{_i}}{G-a} \cup \{a\}$. \inQED\\    

To conclude, let then $S_{x_{_1}}$ be nonelementary and $H_{_1}$ (as in (\sref{1con-1}.A.4.2)) be a trivial $5$-hammock of $G$, by \bfm{(a)} and \bfm{(b)} (above). We may assume that $S_{x_{_2}}$ is nonelementary. If not so, then $J$ is configuration (I) from Figure~\ref{confs} with $x_{_{2}}$ appearing on its outerwalk. Since $C$ is a $4$-face of $G-a$, then by \sref{conf-face}, triviality of $H_{_1}$, (B.2-4), and $4$-connectivity of $G-a$, at least one of $x_{_4},x_{_5}$ also lies on the outerwalk of $J$ and consequently (by planarity of $G-a$) has degree at most $3$ in $G-a$; which is a contradiction. 

By (\sref{1con-1}.A.4.2) and \bfm{(b)} (above), $S_{x_{_2}}$ being nonelementary, implies that there is an $H_{_2}$ (as in (\sref{1con-1}.A.4.2)) that is a trivial $5$-hammock of $G$. Since $C$ is a face, $J$ is configuration (I)-(V). If $J$ is configuration (I), then at least one of $x_{_4},x_{_5}$ lies on the outerwalk of $J$, by \sref{conf-face}. Triviality of $H_{_1},H_{_2}$, and since $a$ is not adjacent to $\{x_{_4},x_{_5}\}$, imply that such a vertex has degree $\leq 4$ in $G$; a contradiction. 

Suppose then that $J$ is configuration (II)-(V). Let $S_{x_{_1}} = \{s,s',t,t'\}, S_{x_{_2}} = \{p,p',q,q'\}$ such that $\{s,s',t,t'p,p',q,q'\}$ are as in \S\ref{salami-and-pie} and Figure~\ref{confs}. Triviality of $H_{_1},H_{_2}$, observation \sref{conf-face}, the fact that $N_G(x_{_j}) \subseteq V(J)$, for $j=4,5$, and since $N_{G}(x_{_i}) = S_{x_{_i}} \cup \{a\}$, $i=1,2$, imply that $J$ has a $3$-disconnector of the form $\{s,p,v\}$ or $\{t,q,v\}$, where $v \in \{x_{_4},x_{_5}\}$, contradicting at least one of (B.2-4).\QED\\

\noindent
\scaps{Proofs.} To conclude our proof of (\sref{1con-1}.A), it remains to prove claims (\sref{1con-1}.A.1-3). This is done next. 

\noindent
\bfm{Proof of} (\sref{1con-1}.A.1)\bfm{.}
By the assumption that $\kappa(G-P) = 1$ and (\ref{setting1}), $G-P$ is the union of an $x_{_4} a$-subchain $X_{_1}$ and an $x_{_5} a$-subchain $X_{_2}$ such that 
$V(X_{_1}) \cap V(X_{_2}) = \{a\}$. Let $y_{_i} \in N_G(a) \cap V(X_{_i})$, $i=1,2$. 
Let $P_{_1}$ be an $x_{_4} y_{_1}$-path in $X_{_1}$ and let $P_{_2}$ be an $x_{_5} y_{_2}$-path in $X_{_2}$. Clearly, $P_{_1}$, and $P_{_2}$ are disjoint and $V(P) \cap V(P_{_i}) = \{x_{_i}\}$, $i=1,2$. 

Assume, towards contradiction, that $(G-a,\{x_{_1},x_{_2},y_{_1},y_{_2}\})$ is a pie and let $F$ be its mold. Since $C$ is a $4$-face and $\kappa(G-a) \geq 4$, we may assume that $\{S_{x_{_1}} \cap F, S_{x_{_2}} \cap F, S_{y_{_1}} \cap F, S_{y_{_2}} \cap F\}$ appear in this clockwise order along $F$. Clearly, $S_{x_{_i}}$, $i=1,2$, if nonelementary, contains $x_{_4}$ and $x_{_5}$ and these appear as specified in (\sref{conf-face}.B). Thus, since $x_{_4},x_{_5} \notin V(P)$, we have that $(V(P) \cap V(P_{_i})) \sm \{x_{_i}\} \not= \emptyset$ for at least one $i =1,2$; contradiction.        
\QED\\

For the proofs of (\sref{1con-1}.A.2-3), we require some preparations. 
Throughout the reminder of this section, we refer to the notation and terminology of \scaps{Setting A} above.
 
Any $4$-disconnector separating $x_{_1}$ from $x_{_2}$ contains $x_{_4}$ and $x_{_5}$; hence, $x_{_4},x_{_5} \in A$. Consequently, $S_{_x}$, if nonelementary, satisfies $\{x_{_4},x_{_5}\} \subset S_{_x} \cap A$. This cannot occur in configuration (I) without contradicting (B.2). 
On the other hand, if $S_{_x}$ is elementary and (thus) $J$ is configuration (I), then $x$ is on the outerwalk of $J$ and is adjacent to two members of $A$. This also cannot occur in configuration (I) without contradicting (B.2). We infer that:
\begin{equation}\label{notI}
\mbox{$J$ is not configuration (I).}
\end{equation}
Since $x_{_1},x_{_2}$ are on the outer face of $G$ (recall $\Pi$), we may assume that $a_{_1} =x_{_4}=s, a_{_4} = x_{_5}=t$, $S_{x} = \{s,s',t,t'\}$, and $S_{y}=\{p,p',q,q'\}$,
where $\{s,s',t,t',p,p',q,q'\}$ are as in \S\ref{salami-and-pie} and Figure~\ref{confs}. 
Let $H_x$ be the $4$-hammock of $J$ satisfying $x \in int_J H_x$ and $\bnd{H_x}{J} = S_x$. In a similar manner define $H_y$ for $y$. Next we prove claims (\sref{1con-1}.A.5-7).\\ 

\noindent
(\sref{1con-1}.A.5) \emph{$V(H_x) = S_x \cup \{x\}$ so that $N_{G-a}(x)=S_x$.}

\proof
By the minimality of $J$ (see (\ref{choice})) and \sref{actually}, $N_G(a) \cap int H_x = \{x\}$. 
Therefore, $\bnd{H_x}{G} = S_x \cup \{x\}$. As we may assume that $H_x$ is not a planar hammock of $G$, it follows that $|V(H_x)| \leq 6$. If equality holds, then there exists a 
(single) vertex $v \in int_G{H_x}$ with $N_G(y) = \bnd{H_x}{G} = S_x \cup \{x\}$; and then $\{x,v,x_{_4},x_{_5}\}$ induce a $K^-_{_4}$. \QED \\

\noindent
(\sref{1con-1}.A.6) \emph{$S_x \cap N_G(a) = \emptyset$.}

\proof
As $\{x_{_4},x_{_5}\} \subset S_{_x} \cap bnd_{G-a} J$, suffices that we show that $s',t' \notin N_G(a)$. Indeed, assume $s' \in N_G(a)$. Then 
$(G-a,\{x,s',\{x_{_1},x_{_2}\} \sm \{x\}, \{\ell,\ell'\} \sm \{y\} \})$ is a Salami with $\{x,s'\}$ at a common end $J'$, by \sref{actually}. The $4$-hammock $J'$ is not configuration (I), by (\ref{notI}). Thus, $J'$ has a $4$-disconnector separating $x$ from $s'$, which is a contradiction as $xs' \in E(G)$, by (\sref{1con-1}.A.5). \QED \\

A consequence of (\sref{1con-1}.A.5-6) is that
\begin{equation}\label{deg}
d_J(s'),d_J(t') \geq 5.
\end{equation}

\noindent
(\sref{1con-1}.A.7) \emph{There are disjoint paths $P'$ and $Q'$ satisfying the following:\\
\indent \emph{(\sref{1con-1}.A.7.1)} $P'$ is an $s'p'$-path and $Q'$ is a $t'q'$-path;\\
\indent \emph{(\sref{1con-1}.A.7.2)} neither meets $A \cup \{p,q,x\} \cup (V(H_{_y}) \sm \{p',q'\})$ (so that $V(P') \cup V(Q') \subseteq J$).\\ 
\indent \emph{(\sref{1con-1}.A.7.3)} Moreover, if $P'$ and $Q'$ are not part of the boundary of a face of $J$ containing \indent \indent \indent \indent \ $\{s',t',p',q'\}$, then at least one of them consists of a single vertex.}\\

\proof
Since $J$ is not configuration (I), by (\ref{notI}), we have that either $\{s',t'\} \cap \{p',q'\} = \emptyset$, and then $\{s',t',p',q'\}$ are cofacial; or $\{s',t'\} \cap \{p',q'\} \not= \emptyset$. We consider these two cases. 

Suppose that there is a face $f$ containing $\{s',t',p',q'\}$. If $V(f) \sm (V(H_{_y}) \cup V(H_{_x}))$ meets $A \cup \{p,q\}$, then the degree of $s'$ or $t'$ is $< 5$; or there is a disconnector $T$, $|T| \leq 3$ separating $y$ from some member of $A$, or there is a disconnector of $G-a$ of size $\leq 3$. These contradict (\ref{deg}), (B.2), and $4$-connectivity of $G-a$, respectively. The paths required then lie on the boundary of $f$; and (\sref{1con-1}.A.7.1) is satisfied.   

If a face as above does not exist, then $|\{s',t'\} \cap \{p',q'\}|\geq 1$. By (\ref{notI}), we may assume then that $p=a_{_2}$ (the case that $q = a_{_3}$ is symmetrical), and therefore $s' \not= p'$, or $\{x,s',p\}$ separate $y$ from $a_{_1}$ in $J$; contradicting (B.2). Consequently, $t'=q'$, and one of the required paths is a single vertex. A path as $P'$ is then contained in an $(s',\{a_{_1},p=a_{_2},t'=q',p'\}$-$4$-fan which exists in $G-a$, contained in $J$, and does not meet $x$, by (\sref{1con-1}.A.5).\QED \\    

%

We are now ready to prove claims (\sref{1con-1}.A.2-3). 

\noindent
\bfm{Proof of} (\sref{1con-1}.A.2)\bfm{.} 
Let $F$ be a $(y,S_y)$-$4$-fan; clearly, $F \subseteq H_y$. Let $P'$ and $Q'$ be an $s'p'$-path and a $t'q'$-paths, respectively, as in (\sref{1con-1}.A.7). If $a_{_2}$ and $a_{_3}$ coincide with $p'$ and $q'$, respectively, then it is easy to see that $Y$ and $Q$ are contained in $F \cup P' \cup Q' \cup \{xs',xt'\}$, respectively. Therefore, let us assume that $a_{_3} \not= q$ (the argument for $a_{_2}$ is symmetrical) which implies that $p$ lies on the outerwalk of $J$ (by the list of possible configurations). 

Let us show that there is an $r \in N_J(a_{_3}) \sm V(Q')$. 
For suppose that $N_J(a_{_3}) \subseteq V(Q')$. Let, then, $b \in N_J(a_{_3})$ such that $t' Q' b$ is minimal. If $P'$ and $Q'$ are parts of the boundary of a face containing $\{s',t',p',q'\}$, then $\{t',b,a_{_4}\}$ is a $3$-disconnector of $J$ contradicting (B.2). To see this, note that $d_J(t') \geq 5$, by (\ref{deg}), and that $N_{G-a}(x) = S_{_x}$. Minimality of $b$ then implies the assertion. 

Consequently, at least one of $P'$ and $Q'$ is a single vertex. Since $p$ is on the outerwalk of $J$, $s'=p'$ implies that $\{x,s'=p',p\}$ is a $3$-disconnector of $J$ contradicting (B.2). 
Thus, $t'=q'=b$ and then $\{x,t'=b=q',a_{_3}\}$ is a $3$-disconnector of $J$ contradicting (B.2).

We have shown that $r$ as above exists. Since $p$ is on the outerwalk of $J$, then $Q' \cup \{x,a_{_4},a_{_3},p\}$ separates $r$ from $a_{_4}$. Also, observe that $r \not=x$, by (\sref{1con-1}.A.5), and that $r \notin int_J H_{_y}$. 

Let $F'$ be as follows: if $q$ is not on the outerwalk of $J$, then $F'$ is an $(r,\{q',q,a_{_4},a_{_2}\})$-$4$-fan. Otherwise, $F'$ is an $(r,\{q',q,a_{_3},a_{_2}\})$-$4$-fan. Let $Q_q$ and $Q_{q'}$ be the members of $F'$ ending at $q$ and $q'$, respectively.\\
(i) By planarity and disjointness of $Q_q$ and $Q_{q'}$, $V(Q_{q'}) \cap V(Q') = \emptyset$.\\
(ii) By definition of $F'$ and planarity, $Q_{q'} \subseteq J -int_J H_y$.

To state $Y$ and $Q$ properly, define $r$ to be $a_{_3}$ if $a_{_3} = q$, and  $r \in N_J(a_{_3}) \sm V(Q')$, otherwise. In an analogous manner define $r'$ for $a_{_2}$ with $p$ replacing $q$. This way we may define a path $P_{p}$ as a $pa_{_2}$-path disjoint of $P'$ and contained in $J-int_J H_{_y}$ in a symmetrical construction to that yielding $Q_{q}$ (above). Planarity asserts that we may choose $P_{p}$ and $Q_{q}$ disjoint. 
Clearly, $Y,Q$ as required are contained in $F \cup P' \cup Q' \cup P_{p} \cup Q_{q} \cup\{xs',xt',ra_{_3},r'a_{_2}\}$.\QED \\

\noindent
\bfm{Proof of} (\sref{1con-1}.A.3)\bfm{.} Let $P',Q'$ be as in the proof of (\sref{1con-1}.A.2). Constructing an $xa_{_2}$-path (resp., $xa_{_3}$-path) containing $P'$ (resp., $Q')$ and $y$ can be done precisely as in the proof of (\sref{1con-1}.2).
We concentrate on constructing an $xa_{_3}$-path not meeting $y$ so that together these paths would constitute the required fan.\\
{\small \scaps{Remark.} In the sole case that $q= a_{_3}$ and $p$ not on the outerwalk of $J$ we shall prefer to construct the paths so that the $xa_{_3}$-path meets $y$; the construction of the other path in this case is symmetrical to what follows.}
  
Let $r \in N_J(a_{_3}) \sm int_J H_y$; such a vertex clearly exists if $a_{_3} \not= q$.
Nonexistence of such a vertex in case $a_{_3} = q$, implies that $\{p,p',q'\}$ is a $3$-disconnector of $J$ separating $y$ from $\{a_{_1},a_{_4}\}$, contradicting (B.2). 

Let $F'$ be an $(r,\{q',q,a_{_3},a_{_4}\})$-$4$-fan. By planarity and the disjointness of the members of $F'$, the path $P_{q'} \in F'$ ending at $q'$ is contained in $J-int_J H_y$, and the required remaining $xa_{_3}$-path is contained in $Q' \cup P_{q'} \cup \{xt',ra_{_3}\}$.\QED

\subsection{Proof of \sref{2con-2}}\label{proof-2con2} 
The reader should be reminded of the \scaps{agreement} specified in \S\ref{pre}.
Throughout this section, $G$ is a $5$-connected nonplanar apex graph satisfying $TK_{_5} \not\subseteq G$; $P$ is an induced $x_{_1} x_{_2}$-path satisfying (\ref{P-disconnected}), where $\{x_{_i}\}^k_{i=1}$ and $K$ are as in (\ref{K}); $H$ is as in (\ref{P-disconnected}); and $a$ is an apex vertex of $G$ as in \sref{2con-2}.

Claim \sref{2con-2} follows from the following two claims; of which (\sref{2con-2}.A) is an exercise.\\
(\sref{2con-2}.A) \emph{$a \in \{x_{_4},x_{_5}\}$; otherwise $TK_{_5} \subseteq G$ or $G-a$ is nonplanar.}\\
(\sref{2con-2}.B) \emph{If $a \in \{x_{_4},x_{_5}\}$, then $TK_{_5} \subseteq G$ or $G$ contains a planar hammock or $K^-_{_4} \subseteq G$.}\\ 

\noindent
\scaps{Common theme.} The following is a common theme of the proofs of (\sref{2con-2}.A-B):
Let $x \in int P$ and let $y,z \in N_H(x)$. 
A $\{yz,x_{_4},x_{_5}\}$-circuit in $H+yz$ clearly implies an $\{x,x_{_4},x_{_5}\}$-circuit in $H+x$ satisfying (F.2). Let then $\X$ be a $\{yz,x_{_4},x_{_5}\}$-frame in $H+yz$, which exists by \sref{frame-separation}. Let $H_{_{y,z}},H_{_{x_{_4}}},H_{_{x_{_5}}}$ be the seeded $\X$-hammocks containing  $\{y,z\},x_{_4},x_{_5}$, respectively. Let $F_{_1}$ and $F_{_2}$ denote the components of such a frame.\\

\noindent
\bfm{Proof of}(\sref{2con-2}.A)\bfm{.}
It is not hard to see that $a \in V(P) \cup \{x_{_3}\}$ contradicts planarity of $G-a$ or $5$-connectivity of $G$. Suppose then that $a \in V(H)$. 
Choose $x,y,z$ as in the \scaps{common theme}; below we refer to the notation and terminology set in the \scaps{common theme}. In $H-yz$ there is an $\{x_{_4},x_{_5}\}$-circuit meeting the ends of $H_{_{y,z}}$ and not meeting $int_{H} H_{_{y,z}}$. Consequently, $H-a-yz$ contains paths $P'$ and $Q$ satisfying:\\
(i) $P'$ is an $x_{_4} x_{_5}$-path,\\
(ii) $Q$ is a $ww'$-path such that $w \in N_{H_{_{y,z}}}(x)$, $w' \in int P'$,\\
(iii) $V(Q) \cap V(P') = \{w'\}$. 
Thus, unless $a \in \{x_{_4},x_{_5}\}$, then a $TK_{_{3,3}}$ is contained in the union of: $P'$, $Q$, $P$, and $V(K) \setminus \{x_{_3}\}$. 
\QED \\ 

\noindent
\bfm{Proof of} (\sref{2con-2}.B) \bfm{.}
Assume, without loss of generality, that $a = x_{_5}$, and choose $x,y,z$ as in the \scaps{common theme}; below we refer to the notation and terminology set in the \scaps{common theme}. We argue as follows: 
\bfm{(a)} $x_{_4} \in int_H H_{x_{_4}}$, by definition of frames, so that $|V(H_{x_{_4}})| \geq 3$. \\
%
%
\bfm{(b)} By $5$-connectivity, $H_{x_{_4}}$, having order $\geq 3$, has at least two attachments on $P$. Let $Q$ denote the extremal segment of $H_{x_{_4}}$ on $P$;
such does not consist of a single edge for then the ends of $Q$ and the ends of $H_{x_{_4}}$  form a $4$-disconnector of $G$.\\ 
\bfm{(c)} $N_G(x_{_5}) \cap int Q = \emptyset$. For suppose $\ell \in N_G(x_{_5}) \cap int Q$. Then, due to planarity, there exists an $\ell' \in N_{H_{x_{_4}}}(\ell)$. An $\ell' x_{_5}$-path containing $x_{_4}$ clearly exists in $H$. Indeed, a seeded-$\ell'$-path in $H_{_{x_{_4}}}$ and a seeded-through path in $H_{x_{_5}}$ are easily extended into such a path. Consequently, an $\{\ell,x_{_4},x_{_5}\}$-circuit satisfying (F.2) exists in $H+\ell$.\\ 
\bfm{(d)} Since the ends of $Q$ and the ends of $H_{x_{_4}}$ do not form a $4$-disconnector of $G$, $x_{_3}$ is adjacent to $int Q$. We may assume that $x_{_3}$, the ends of $Q$, and the ends of $H_{x_{_4}}$ do not form the boundary of a planar hammock.
Since $x_{_4} \in int_H H_{x_{_4}}$, then these form the boundary of a $5$-hammock of $G$ of order $6$ with $x_{_4}$ as the sole vertex in the interior of this hammock. This implies that $x_{_4} x_{_3} \in E(G)$; contradicting (\ref{stable}).\QED

\subsection{Proof of \sref{1con-2}} Throughout this section, $G$ is a $5$-connected nonplanar apex graph satisfying $TK_{_5} \not\subseteq G$; $P$ is an induced $x_{_1} x_{_2}$-path satisfying (\ref{P-disconnected}), where $\{x_{_i}\}^k_{i=1}$ and $K$ are as in (\ref{K}); $H$ is as in (\ref{P-disconnected}); and $a$ is an apex vertex of $G$ as in \sref{1con-2}. 

Claim \sref{1con-2} follows from the claims (\sref{1con-2}.A-B). Tough not necessarily short, (\sref{1con-2}.A) is an exercise which we omit. \\
(\sref{1con-2}.A) \emph{If $a \notin \{x_{_4},x_{_5}\}$, then $G-a$ is nonplanar.}\\
(\sref{1con-2}.B) \emph{If $a \in \{x_{_4},x_{_5}\}$, then $TK_{_5} \subseteq G$ or $K^-_{_4} \subseteq G$ or $G$ has a planar hammock.}\\ 

\noindent
\bfm{Proof of} (\sref{1con-2}.B) \bfm{.}
Assume, without loss of generality, that $a = x_{_5}$. The circuit $C = P+x_{_1}x_{_4}+x_{_2}x_{_4}$ separates $x_{_3}$ from $H-a$. Fix an embedding of $G-a$ such that $x_{_3}$ lies in $ext C$ and $H-a$ lies in $int C$. Let $B$ be the leaf block  
of $H$ containing $x_{_4}$. Let $H'$ be the subchain of $H$ containing all blocks of $H$ other than $B$. Let $y$ denote the cut vertex of $H$ present in $B$; by \sref{2con-2}, $H'$ and $y$ exist.\\

\noindent
\scaps{Outline.} We now outline our argument for proving \sref{1con-2}; adjourning technical details to the end of the section. Proofs of claims hereafter follow below. Initially we prove that\\

\noindent
(\sref{1con-2}.B.1) \emph{$B$ is nontrivial (i.e., does not consist of a single edge).}\\

\noindent
By (\sref{1con-2}.B.1), the boundary of the outer face of $B$ is a circuit comprised of two 
$x_{_4}y$-paths. An attachment vertex of $B$ in $int P$ is incident to at least one of these paths. Planarity asserts that there is at most one vertex in $int P$ that is adjacent to both of these paths. Consequently, there are attachments $z,z' \in int P$ of $B$ such that $zPx_{_2}$ contains all attachments adjacent to one of these $x_{_4} y$-paths, and $z'Px_{_1}$ contains all attachments of $B$ adjacent to the other path. These segments of $P$ are disjoint, unless $z = z'$. We will prove that\\

\noindent  
(\sref{1con-2}.B.2) $N_G(x_{_5}) \cap (z'Px_{_1}) = N_G(x_{_5}) \cap (zPx_{_2}) = \emptyset$.\\
(\sref{1con-2}.B.3) $|\{z,z',x_{_1},x_{_2}\}|=4$.\\
(\sref{1con-2}.B.4) \emph{$N_G(x_{_3}) \cap (z'Px_{_1}) \not= \emptyset$ or $ N_G(x_{_3}) \cap (zPx_{_2}) \not= \emptyset$.}\\
(\sref{1con-2}.B.5) \emph{$H'$ does not consist of a single edge (i.e. $x_{_5}y$).}\\

Let $w,w' \in V(P)$ be attachments of $H'-x_{_5}$ on $P$ such that $wPx_{_2}$ is minimal and 
$w'Px_{_1}$ is minimal; these exist by (\sref{1con-2}.B.5). By planarity of $G-a$, $w,w' \in zPz'$. We will show that\\

\noindent
(\sref{1con-2}.B.6) \emph{$N_G(x_{_3}) \cap (wPw') \not= \emptyset$.}\\

Let $h_{_1} \in N_G(x_{_3}) \cap (wPw')$, by (\sref{1con-2}.B.6). 
Let $h'_{_1} \in N_G(h_{_1}) \cap V(H'-x_{_5})$. 
Let $h_{_2} \in N_G(x_{_3}) \cap (z'Px_{_1})$, without loss of generality, by (\sref{1con-2}.B.4). 
By (\sref{1con-2}.B.2), $h_{_2}$ is at least $5$-valent in $G-x_{_5}$. Planarity then implies that there is an $h'_{_2} \in N_G(h_{_2}) \cap V(B) \sm \{x_{_4},y\}$. Indeed, if $h_{_2} x_{_4} \in E(G)$, then since $P$ is induced $x_{_1}$ is at most $4$-valent in $G$; contradicting $5$-connectivity. Let $b \in N_G(x_{_2}) \cap V(B) \sm \{y,x_{_4}\}$. 
Such a vertex exists for otherwise $x_{_2}y \in E(G)$, implying that $z =x_{_2}$; contradicting 
(\sref{1con-2}.B.3). We will see that\\

\noindent
(\sref{1con-2}.B.7) \emph{$B-x_{_4}-y$ contains a $bh'_{_2}$-path.}\\

With the above claims in place, we construct a $TK_{_5}$ in $G$ whose branch vertices are $x_{_3},h_{_1},h_{_2},x_{_1},x_{_2}$ as follows. Since $\{h_{_1},h_{_2},x_{_1},x_{_2}\} \subset N_G(x_{_3})$, suffices that we show that $G-x_{_3}$ has a $TK_{_4}$ rooted at $\{h_{_1},h_{_2},x_{_1},x_{_2}\}$. Since $C$ contains all these vertices, such a rooted 
$TK_{_4}$ is implied by existence of an $x_{_1} h_{_1}$-path $Q_{_1}$ and an $x_{_2}h_{_2}$-path $Q_{_1}$ in $G-x_{_3}$ that are internally-disjoint of $C$ and disjoint  from one another. 
Put $Q_{_1} = x_{_1}x_{_5} + x_{_5}h_{_3} + Q'_{_1} +h_{_1}h'_{_1}$, where $h_{_3} \in N_{H'}(x_{_5})$ and $Q'_{_1}$ is an $h'_{_1}h_{_3}$-path in $H'-x_{_5}$. 
The path $Q'_{_1}$ exists since is not a cut vertex of $H'$ as it is contained internally in one of its leaf blocks. Next, put $Q_{_2} = x_{_2}b + Q'_{_2} + h'_{_2}h_{_2}$, where $Q'_{_2}$ is as in (\sref{1con-2}.B.7).\\ 

\noindent
\scaps{Proofs.} It remains to prove (\sref{1con-2}.B.1-7). 

\noindent
\bfm{Proof of} (\sref{1con-2}.B.1)\bfm{.} 
If $B$ is trivial, then it consists of the edge $x_{_4} y$, as $x_{_4}$ is contained internally in $B$. Planarity, $P$ being induced, the $P$-bridge meeting $x_{_3}$ is a singleton, imply that there is a $w \in \{x_{_1},x_{_2}\}$ and a vertex in $int P$ such that these two vertices and $x_{_4}$ define a triangle. Planarity then implies that $wy \notin E(G)$ and thus $w$ is at most $4$-valent in $G$; contradicting $5$-connectivity.\inQED \\

\noindent
\bfm{Proof of} (\sref{1con-2}.B.2)\bfm{.}
Suppose, to the contrary, that $r \in N_G(x_{_5}) \cap (zPx_{_2})$. Clearly, $N_B(r) \not= \emptyset$. Suffices that we prove that 
$N_B(r) \sm \{y\} \not= \emptyset$. Indeed, if so, and $r' \in  N_B(r) \sm \{y\}$, then 
an $(x_{_4},\{y,r'\})$-$2$-fan in $B$ (which exists as $B$ is nontrivial by (\sref{1con-2}.A)), a $y \ell$-path in $H'-x_{_5}$, where $\ell \in N_{H'-x_{_5}}(x_{_5})$,
and the edges $\{x_{_5} r, rr', x_{_5} \ell\}$ define a circuit satisfying (F.2).

Suppose then that $N_B(r) = \{y\}$. If so, then the following terms are satisfied:\\
(i) $\{y,z,r\}$ is a triangle; indeed, $rz \notin E(P)$ implies that $zPr \subset N_G(y)$ and thus that $K^-_{_4} \subset G$.\\
(ii) $rx_{_3} \in E(G)$; implied by planarity and the assumption that $N_B(r) =\{y\}$.\\
(iii) $zx_{_3} \notin E(G)$; or $K^-_{_4} \subset G$ by (i) and (ii) above.\\
(iv) $z \not= z'$; otherwise $H'$ consists of the single edge $x_{_5} y$ and thus 
$\{z,y,r,x_{_5}\}$ induce a $K^-_{_4}$.\\
(v) $N_G(x_{_3}) \cap (zPx_{_1}) = \emptyset$, meaning that $N_G(x_{_3}) = x_{_1} \cup Y$ such that $Y \subseteq rPx_{_2}$. To see this, we show that if there is a $u \in N_G(x_{_3}) \cap (zPx_{_1})$, then $TK_{_5} \subseteq G$ which is constructed as follows:

(v.a) $\{r,x_{_5},x_{_2},x_{_3}\}$ are the branch vertices of a $TK_{_4}$. Indeed,  $\{r,x_{_5},x_{_2},x_{_3}\}$ induce a $C_{_4}$ \indent \indent \ \  and the paths $rPx_{_2}$ and $(x_{_5},x_{_1},x_{_3})$ complete this $C_{_4}$ into a $TK_{_4}$.

(v.b) Let $F =\{P_{_1},P_{_2},P_{_3},P_{_4}\}$ be the following $(y,\{r,x_{_5},x_{_2},x_{_3}\})$-$4$-fan: $P_{_1}$ consists of the \indent \indent \ \ edges 
$\{yz,rz\}$. $P_{_2}$ is a $y x_{_5}$-path in $H'$. Let $P'_{_3}$ and $P'_{_4}$ be the 
the members of a $2$-fan \indent \indent \ \ in $B$ from $y$ to $\{x_{_4},\ell\}$, where 
$\ell \in N_B(z') \sm \{x_{_4}\}$; such a vertex exists for otherwise \indent \indent \ \ planarity would imply that $x_{_1}$ is at most $4$-valent in $G$.

(v.c) Let $P_{_3}$ and $P_{_4}$ be obtained from $P'_{_3}$ and $P'_{_4}$ by adding the edge 
$x_{_4}x_{_2}$ to the member \indent \indent \ \ of $\{P'_{_3},P'_{_4}\}$ that is a $yx_{_4}$-path. The other member of $\{P'_{_3},P'_{_4}\}$ is a $y \ell$-path. Add to it \indent \indent \ \ the edge $z' \ell$ and the subpath $z'Pu$ and the edge $x_{_3}u$.

\noindent
The union of $F$ and the $TK_{_4}$ constructed in (v.a) form a $TK_{_5}$.

We may now assume that $N_G(x_{_5}) \cap (z' P x_{_1}) = \emptyset$; otherwise, by a symmetrical argument to the one above we will have that $N_G(x_{_3}) \cap (z'Px_{_2}) = \emptyset$
contradicting (ii). Thus, $\{z',y,r\}$ is a $3$-disconnector of $G-a$ separating $x_{_3}$ from $H'-x_{_5}$ which is nonempty as seen in (iv); contradiction. \inQED \\

\noindent
\bfm{Proof of} (\sref{1con-2}.B.3)\bfm{.} If $z =z'$, then $H'$ consists of the single edge $x_{_5} y$; consequently, $x_{_5}$ is adjacent to at least one of $(z'Px_{_1})$ or $(zPx_{_2})$
contradicting (\sref{1con-2}.B.2). Next, if, say, $x_{_1} = z'$, then either $\{x_{_1},y,z,x_{_2}\}$ is a $4$-disconnector of $G$ or  $\{x_{_1},y,z,x_{_2},x_{_3}\}$ is the boundary of a planar hammock since $B$ is nontrivial, by (\sref{1con-2}.B.1).\inQED \\

\noindent
\bfm{Proof of} (\sref{1con-2}.B.4)\bfm{.} Suppose not. Then $\{x_{_1},z,y,z',x_{_2}\}$ are the boundary of a $5$-hammock of $G$ of order $6$ (otherwise we have a planar hammock of $G$);
implying that $B$ is trivial, contradicting (\sref{1con-2}.B.1).\inQED \\

\noindent
\bfm{Proof of} (\sref{1con-2}.B.5)\bfm{.} Suppose $H'$ consist of the single edge $x_{_5}y$.
The segment $z'Pz$ consists of a the single edge $zz'$ then; for vertices in $(z'Pz)$ can only have degree $\leq 4$ since $x_{_3}$ is a singleton $P$-bridge. Thus, $z,z' \in N_{G-x_{_5}}(x_{_5})$, by (\sref{1con-2}.B.2); implying that $\{z,z',x_{_5}\}$ is a triangle. 
Thus, neither of $z,z'$ is adjacent to $y$; or $K^-_{_4} \subseteq G$. In addition, at least one of $z,z'$, say $z$, is not adjacent to $x_{_3}$; or $K^-_{_4} \subseteq G$. To be at least $5$-valent, then, $|N_B(z)| \geq 2$. 
We then construct a $\{z,x_{_5},x_{_4}\}$-circuit satisfying (F.2) as follows. The union of
the edges $\{x_{_5}y,x_{_5}z,z\ell\}$, where $\ell \in N_B(z) \sm \{y\}$,
and an $(x_{_4},\{y,\ell\})$-$2$-fan in $B$ form such a circuit. \inQED \\

\noindent
\bfm{Proof of} (\sref{1con-2}.B.6)\bfm{.} Suppose not. The set $\{x_{_5},y,w,w'\}$ has size $\leq 4$. Since this set is not a $4$-disconnector, $H'$ consists of the single edge $x_{_5}y$; contradiction to (\sref{1con-2}.B.5). \inQED \\ 


\noindent
\bfm{Proof of} (\sref{1con-2}.B.7)\bfm{.} 
If such a path does not exist, then $\{x_{_4},y\}$ is a $2$-disconnector of $B$ separating $h'_{_2}$ from $b$. That is, $B = B_{_1} \cup B_{_2}$ such that $G[\{x_{_4},y\}] = B_{_1} \cap B_{_2}$, $h'_{_2} \in V(B_{_1})$ and $b \in V(B_{_2})$. Since $|\{x_{_1},x_{_2},z,z'\}| =4$, by (\sref{1con-2}.B.3), $N_{B_{_i}}(x_{_i}) \sm \{y,x_{_4}\} \not= \emptyset$ for $i=1,2$. 
As neither of $\{x_{_1},z',y,x_{_4}\}$, $\{x_{_2},z,y,x_{_4}\}$ is a $4$-disconnector of $G$,
it follows that $x_{_3}$ is adjacent to both $(x_{_1}Pz')$ and $(x_{_2}Pz)$. 
Since we assume that neither of $\{x_{_3},x_{_1},z',y,x_{_4}\}$, $\{x_{_3},x_{_2},z,y,x_{_4}\}$
form the boundary of a planar hammock, each of these sets is the boundary of a $5$-hammock of $G$ of order $6$. This then implies that $x_{_3}$ is adjacent to $B$; contradiction to the premise of the \scaps{Disconnected Case}.\inQED \\
\QED


\begin{thebibliography}{99}


\bibitem{diestel}
R. Diestel, \emph{Graph Theory}, third edition, Springer, $2005$. 

\bibitem{MaYu1}
J. Ma and X. Yu, Independent paths and $K_{_5}$-subdivisions, J. Combinatorial Theory B (to appear).

\bibitem{MaYu2}
J. Ma and X. Yu, $K_{_5}$-subdivisions in graphs containing $K^-_{_4}$, submitted manuscript.  

\bibitem{Mader1}
W. Mader, $3n-5$ edges do force a subdivision of $K_{_5}$, Combinatorica, $\mbs{18}$ ($4$) ($1998$) $569$-$595$.


\bibitem{mw}
D. M. Mesner and M. E. Watkins, Cycles and connectivity in graphs.
Can. J. Math.,  $\mbs{19}$ ($1967$) $1319-1328$.

\bibitem{mt}
B. Mohar and C. Thomassen, \emph{Graphs on Surfaces}, $2001$, The Johns Hopkins University Press Baltimore and London.

\bibitem{seymour}
P. D. Seymour, Disjoint paths in graphs, Discrete Math., $\mbs{29}$ ($1980$) $293-309$. 

\bibitem{thomassen}
C. Thomassen, $2$-linked graphs, European J. Combin., $\mbs{1}$ ($1980$) $371-378$.

\bibitem{Yu}
X. Yu, Subdivisions in planar graphs, J. Combin. Theory Ser. B, $\mbs{72}$ ($1998$) $10-52$. 

\end{thebibliography}
\end{document}